\documentclass[journal,twoside,web]{ieeecolor}

\usepackage{generic}
\usepackage{cite}
\usepackage{amsmath,amssymb,amsfonts}
\usepackage{graphicx}
\usepackage{textcomp}
\usepackage{stmaryrd}
\usepackage{tabularx}

\usepackage[font=footnotesize]{caption}
\usepackage{subcaption}
\usepackage{color}

\newcommand{\f}[2]{\frac{#1}{#2}}

\newcommand{\mbf}[1]{\mathbf{ #1}}
\newcommand{\tbf}[1]{\textbf{#1}}

\newcommand{\mcl}[1]{\mathcal{#1}}
\newcommand{\R}{\mathbb{R}}
\newcommand{\N}{\mathbb{N}}
\renewcommand{\S}{\mathbb{S}}

\newcommand{\bmat}[1]{\begin{bmatrix}#1\end{bmatrix}}

\newcommand{\smallbmat}[1]{\left[\scriptsize\begin{smallmatrix}
#1\end{smallmatrix} \right]}
\newcommand{\mat}[1]{\begin{matrix}#1\end{matrix}}

\newcommand{\semismallbmat}[1]{\footnotesize\bmat{#1}}

\let\bl\bigl
\let\bbl\Bigl

\let\bbbbl\Biggl
\let\br\bigr
\let\bbr\Bigr

\let\bbbbr\Biggr

\title{\LARGE \bf
Efficient Data Structures for 
 Representation of Polynomial Optimization Problems: Implementation in SOSTOOLS
}

\author{Declan Jagt, Sachin Shivakumar, Peter Seiler, Matthew Peet %
\thanks{\tbf{Acknowledgement:} This work was supported by National Science Foundation grants CMMI-1935453 and CMMI-1931270.} %
}

\begin{document}

\maketitle
\thispagestyle{empty}
\pagestyle{empty}

\begin{abstract}

We present a new data structure for representation of polynomial variables in the parsing of sum-of-squares (SOS) programs. In SOS programs, the variables $s(x;P)$ are polynomial in the independent variables $x$, but linear in the decision variables $P$. Current SOS parsers, however, fail to exploit the semi-linear structure of the polynomial variables, treating the decision variables as independent variables in their representation. This results in unnecessary overhead in storage and manipulation of the polynomial variables, prohibiting the parser from addressing larger-scale optimization problems.
To eliminate this computational overhead, we introduce a new representation of polynomial variables, the ``dpvar'' structure, that is affine in the decision variables. We show that the complexity of operations on variables in the dpvar representation scales favorably with the number of decision variables. We further show that the required memory for storing polynomial variables is relatively small using the dpvar structure, particularly when exploiting the MATLAB sparse storage structure. Finally, we incorporate the dpvar data structure into SOSTOOLS 4.00, and test the performance of the parser for several polynomial optimization problems.

\end{abstract}

\section{INTRODUCTION}

Many problems in analysis and control of nonlinear systems can be formulated as polynomial optimization problems.
Since testing nonnegativity of polynomials is NP-hard~\cite{blum1998Complexity}, polynomial constraints of the form $s(x)\geq0$ for all $x\in\R^n$ are often tightened to sum-of-squares (SOS) constraints $s\in\Sigma_s$, where $\Sigma_s$ denotes the set of functions that may be expanded as $s(x)=\sum_{i}^{} p_i(x)^2$ for some polynomial  functions $p_i\in\R[x]$. Feasibility of $s\in\Sigma_s$ in turn is equivalent to existence of a positive semidefinite matrix $Q\geq 0$ and a vector of monomials $Z_{d}$ such that $s(x)=Z_{d}(x)^T Q Z_{d}(x)$, allowing SOS constraints to be expressed as LMIs. In this manner, SOS programs (SOSPs) can be formulated as semidefinite programs (SDPs), which may be solved in polynomial time~\cite{boyd1994LMIs}. For recent applications of SOS programming, see~\cite{liu2020PermissiblePerturbation,li2021TrackingControl,wang2020AttractionDomains}.

\begin{figure}[h!]
	\centering
	\vspace*{-0.3cm}
	\hspace*{-0cm}\includegraphics[width=0.5\textwidth]{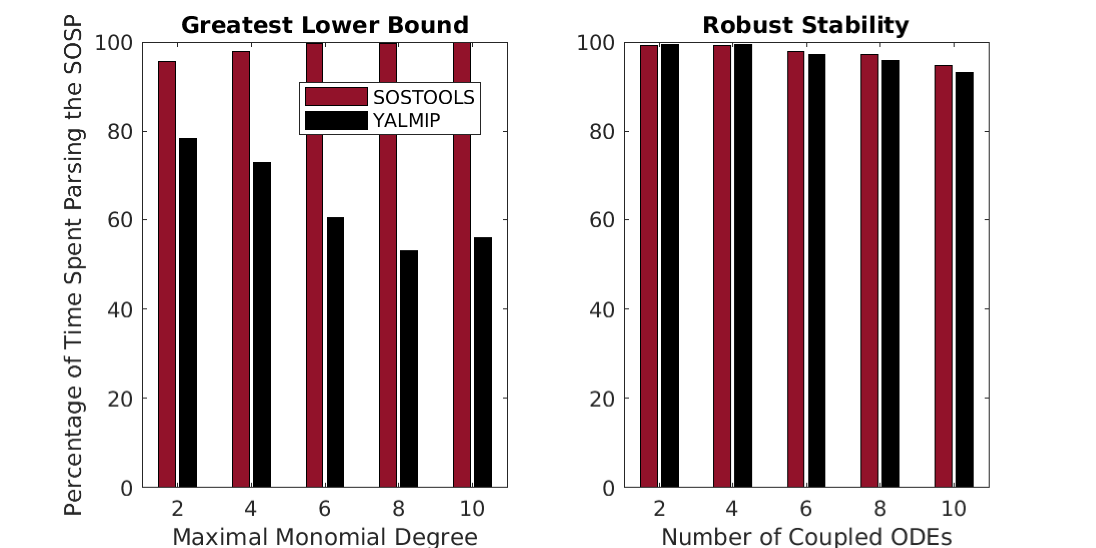}
	\caption{\footnotesize Percentage of execution time spent parsing the greatest lower bound problem from Subsection~\ref{sec:subsec_Greatest_Lower_Bound} (Eqn.~\eqref{eq:GLB_SOS}) and the robust stability problem from Subsection~\ref{sec:subsec_Robust_Stability} (Eqn.~\eqref{eq:robust_stability_SOS}), using SOSTOOLS 3.04 and YALMIP. Using either parser, less than 50\% of the time spent on each problem is actually spent on solving the associated SDP, with the parsing of the robust stability program even taking up more than 90\% of the time.}
	\label{fig:parse_time_ratio_intro}
	\vspace*{-0.3cm}
\end{figure}

The typical process of numerically solving SOSPs consists of two stages: the \textit{parsing} of the SOSP, i.e. the implementation of the program and conversion to an SDP; and the actual \textit{solving} of this SDP. Unfortunately, the computational complexity associated with both of these stages increases rapidly with the size of the SOSP, as a result of which many large-scale applications of SOS programming remain unsolvable. This failure to tackle large-scale problems has prompted several variations on SOS programming to be proposed, reducing complexity of the problem by imposing more restrictive constraints on the positive semidefinite matrix $Q$~\cite{ahmadi2014DSOS,waki2006SSOS,zheng2019SparseDSOS}. However, the goal of these modifications is primarily to reduce the computational complexity of the solving stage of the SOS programming process, offering little to no reduction in the cost of parsing the SOSP. As such, even if larger-scale problems can be solved with these modifications, the computational cost of parsing such programs may still make numerical implementation impossible. In fact, in many cases, the computational complexity of parsing the SOSP far exceeds that associated to solving the resulting SDP (see Fig.~\ref{fig:parse_time_ratio_intro}), a discrepancy that will only be exacerbated by reducing the complexity of the SDP. 


For the greatest lower bound problem and robust stability test presented in Subsection~\ref{sec:subsec_Greatest_Lower_Bound} and~\ref{sec:subsec_Robust_Stability}, Fig.~\ref{fig:parse_time_ratio_intro} shows what percentage of the time required to solve each problem is spent on parsing the SOSP. Results are shown using the well-established SOS parsers SOSTOOLS 3.04~\cite{prajna2002SOSTOOLS} and YALMIP~\cite{lofberg2004YALMIP} to parse the problems, using SEDUMI~\cite{sturm1999SEDUMI} to solve the resulting SDP. 
The results show that both parsers consistently require more time to construct the SDP from the SOSP than it takes to actually solve this SDP, frequently spending more than 90\% of the execution time on parsing. 
In this paper, we show that the percentage of the time spent on parsing can be significantly reduced, proposing a new representation of polynomial variables that allows for more efficient parsing of SOSPs.

In converting an SOSP to an SDP, SOS parsers use finite monomial bases $Z_d$ to represent the polynomial variables. Here, we let $Z_{d}\in\R^{n_1}[x]$ denote a vector containing all monomials in variables $x_1,\hdots,x_p$ of degree at most $d$, where $n_1:=\frac{(p+d)!}{p!d!}$. These monomials may be numerically represented as a matrix $Z_{\text{M},d}\in\N^{n_1\times p}$ containing the degrees of each variable in each monomial, so that e.g.
\vspace*{-0.3cm}
\begin{align*}
Z_2(x_1,x_2)&=\bmat{1 \\x_2\\ x_2^2\\ x_1\\ x_1x_2 \\x_1^2}	&	&\text{and}	&
Z_{\text{M},2}&=\overbrace{\bmat{0&\ 0\\0&\ 1\\0&\ 2\\1&\ 0\\1&\ 1\\2&\ 0}}^{[x_1,x_2]}.
\end{align*}
Using such a monomial basis, an SOS variable $s\in\Sigma_s$ of degree at most $2d$ can be represented in the quadratic form \vspace*{-0.05cm}
\begin{align*}
s(x;Q)&=Z_{d}(x)^T Q Z_{d}(x),
\end{align*}
where now $Q\in\S^{n_1\times n_1}$ is a \textit{decision variable}. Meanwhile, any polynomial $p\in\R[x]$ of degree $2d$ is uniquely defined by a vector of coefficients $c\in\R^{n_2}$ for $n_2:=\frac{(p+2d)!}{p!(2d)!}$, and may be represented in the linear \textit{pvar} form as \vspace*{-0.05cm}
\vspace*{-0.0cm}
\begin{align}\label{eq:pvar_intro}
p(x)&=c^T Z_{2d}(x).
\end{align}
Finally, interface with SDP solvers requires polynomial constraints $g(x;\xi)=0$, parameterized by decision variables $\xi$, to be expressed in the SDP format
\begin{align*}
0=g(x;\xi)&=(A\xi - b)^TZ(x),	& & \text{imposing} 	& A\xi&=b.
\end{align*}
For example, letting $s_1(x_1;\xi)\!=\!\smallbmat{1\\x_1}^T\smallbmat{\xi_1&\xi_2\\\xi_2&\xi_3}\smallbmat{1\\x_1}$ for $\smallbmat{\xi_1&\xi_2\\\xi_2&\xi_3}\geq 0$, 
and defining $p_1(x_1)\!:=\!1-2x_1^2$,
the constraint 
\begin{align*}
0&=g_1(x_1;\xi):=s_1(x_1;\xi)p_1(x_1)-1+4x_1^4,	&
\end{align*}
can be equivalently represented in the SDP format as \vspace*{-0.1cm}
\begin{align*}
0=g_1(x_1;\xi)
&=\bbbbl(\underbrace{\semismallbmat{1&0&0\\0&2&0\\-2&0&1\\0&-4&0\\0&0&-2}}_{A}\underbrace{\semismallbmat{\xi_1\\\xi_2\\\xi_3}}_{\xi}-\underbrace{\semismallbmat{1\\0\\0\\0\\-4}}_{b}\bbbbr)^T\semismallbmat{1\\x_1\\x_1^2\\x_1^3\\x_1^4}.	\\[-2.0em]
\nonumber
\end{align*}
In order to derive this expression, however, an SOS parser would have to compute the product $s_1(x;\xi)p_1(x)$ without knowing the values of the decision variables $\xi$. To this end, the approach of current parsers is to treat the decision variables as independent variables, and represent SOS variables $s$ in the linear form as
\begin{align*}
s(x;\xi)=c^T \bar{Z}_{2d}(x;\xi)
\end{align*}
where $\bar{Z}_{2d}(x;\xi)\!:=\!\semismallbmat{1\\\xi}\otimes Z_{2d}(x)$ is now a vector of monomials in the joint set of variables $(x,\xi)$ -- meaning $Z_{2d}$ will be rather long. Although this linear format allows operations such as multiplication to be performed relatively easily, using e.g.
\begin{align*}
&c_1^T Z_{2d}(x)c_2^T Z_{2d}(x;\xi)=(c_1\otimes c_2)^T\bl(Z_{2d}(x)\otimes Z_{2d}(x;\xi)\br),
\end{align*}
the complexity of operations like multiplication will scale poorly with the number of decision variables $\xi$. Moreover, once the constraint has been converted to one of the form $0=c^T\bar{Z}(x;\xi)$, substantial computational effort may still be required to extract the decision variables $\xi$ from $\bar{Z}$, and define the necessary matrix $A \in \R^{m\times n}$ and vector $b \in \R^{m}$ to express the constraint in the SDP format $0=(A\xi-b)^T Z(x)$.


To reduce the computational overhead associated with parsing SOS programs, we propose a new representation of polynomial decision variables which tracks more closely with the SDP constraint format, while allowing for efficient conceptual and numerical manipulation of the resulting polynomial objects. Specifically, we represent a polynomial variable $s\in\R[x;\xi]$, parameterized by decision variables $\xi$ as
\begin{align}\label{eq:scalar_dpvar}
 s(x;\xi)&:=Z_1(\xi)^T C Z_{d}(x)=\bmat{1\\ \xi}^T C Z_{d}(x),
\end{align}
so that, for example
\vspace*{-0.65cm}
\begin{align*}
s_1(x_1;\xi)&
=\semismallbmat{1\\x_1}^T\semismallbmat{\xi_1&\xi_2\\\xi_2&\xi_3}\semismallbmat{1\\x_1}
={\scriptsize\bmat{1\\\xi_1\\\xi_2\\\xi_3}}^T\overbrace{\scriptsize\bmat{0&0&0\\1&0&0\\0&2&0\\0&0&1}}^{C}\overbrace{\scriptsize\bmat{1\\ x_1\\x_1^2}}^{Z_2(x_1)}.
\end{align*}
We refer to this variable structure as the \textit{decision polynomial variable}, or \textit{dpvar} representation --  a generalization of the linear \textit{polynomial variable}, or \textit{pvar} representation to polynomials with decision variables. As will be shown in Section~\ref{sec:dpvar_operations}, use of this format accounts for linearity with respect to the decision variables and eliminates polynomial manipulations involving decision variables. Furthermore, in this format, translation of an equality constraint such as $s(x;\xi)=0$ to SDP format is trivial, in that
\begin{align*}
s(x;\xi)=\bmat{1\\ \xi}^T C Z_{d}(x)&=\bmat{1\\ \xi}^T \bmat{c_1^T\\C_2^T} Z_{d}(x)\\&=(\xi^TC_2^T+c_1^T)Z_{d}(x),
\end{align*}
so that $s=0$ may be equivalently expressed as an LMI constraint $C_2\xi=-c_1$. Furthermore, by eliminating the need for construction of extremely large transition matrices, memory requirements are significantly reduced. Finally, while the resulting $C$ matrices are still rather large (as is required for densely-defined polynomial expressions), when the number of terms in these matrices is small, the dpvar structure exploits the sparse matrix representation features of MATLAB to dramatically reduce computation time - see Section~\ref{sec:MATLAB_sparsity}.

In the remainder of this paper, we carefully detail and analyze how an ideal parser should integrate the dpvar structure into the parsing of SOS optimization problems. Specifically, an ideal parser should
\begin{enumerate}
  \item Exploit structure in polynomial computations. In particular, for polynomial multiplication, addition, substitution, etc., the parser should exploit the affine appearance of the decision variables to reduce computational overhead. 
  \item Be based on analytic expressions for the mathematical operations.
  \item Allow for fully dense polynomial structures.
  \item Make efficient use of the platform-specific sparsity structure to minimize memory usage and computational complexity for sparse polynomial objects.
  \item Be scalable to hundreds of thousands of decision variables.
\end{enumerate}
In the following sections, we show how the dpvar structure can be used to achieve these goals in the context of the MATLAB programming language and associated sparsity package.

\section{Preliminaries}

\subsection{Notation}

We denote $\R^{m\times n}[x;\xi]$ as the set of $m\times n$ matrix-valued polynomials in variables $x$ and $\xi$. We denote $Z_{d}\in\R^{n}[x]$ as a vector consisting of all monomials in $x$ up to degree $d$, and $\hat{Z}_d\subseteq Z_{d}$ as a vector consisting of only a subset of these monomials. We will often refer to $Z_{d}$ in terms of the degrees of the variables appearing in each monomial, so that e.g.
\begin{align*}
x_1^2 x_2 x_4^4 = \overbrace{\bmat{2&1&0&4}}^{[x_1,x_2,x_3,x_4]}.
\end{align*}
For any monomial basis $Z_d$, we use $Z_{\text{M},d}\in\N^{n\times p}$ to denote the associated matrix of degrees, where $\N$ denotes the set of nonnegative integers and $p$ the number of independent variables. We let $nnz(A)$ denote the number of nonzero elements of a (sparse) matrix $A\in\R^{n\times m}$.  We use big O notation $f(N)=\mcl{O}(g(N))$ for scalar functions $f,g$ to indicate that there exists some constant $C>0$ such that $|f(N)|\leq Cg(N)$ for all $N\in\R$.

\subsection{Example Polynomials}

Throughout the paper, various concepts will be illustrated using the example polynomial $p_1(x_1)=1-2x_1^2$, and the polynomial variable $s_1(x_1;\xi)=\smallbmat{1\\\xi_1}^T\smallbmat{\xi_1&\xi_2\\\xi_2&\xi_3}\smallbmat{1\\x_1}$. Here, the polynomial $p_1$ can be represented in terms of the monomial vector $Z_{2}(x_1)$ in the pvar format as
\begin{align}\label{eq:example_p1_pvar}
p_1(x_1)=b_1^T Z_2(x_1)=\underbrace{\semismallbmat{1\\0\\-2}^T}_{b_1^T}\underbrace{\semismallbmat{1\\x_1\\x_1^2}}_{Z_2(x_1)}.
\end{align}
Similarly, the polynomial variable $s_1$ can be represented in terms of the monomial vectors $Z_1(\xi)$ and $Z_{2}(x_1)$ in the dpvar representation as \vspace*{-0.3cm}
\begin{align}\label{eq:example_s1_dpvar}
s_1(x_1;\xi) &=Z_1(\xi)^T C_1 Z_2(x_1)
=\underbrace{{\footnotesize\bmat{1\\\xi_1\\\xi_2\\\xi_3}}^T}_{Z_1(\xi)^T}\underbrace{\footnotesize\bmat{0&0&0\\1&0&0\\0&2&0\\0&0&1}}_{C_1}{\footnotesize\bmat{1\\ x_1\\x_1^2}}.
\end{align}
or in terms of the monomial vector $\bar{Z}_{2}(x_1;\xi)$ in the pvar representation as \vspace*{-0.3cm}
\begin{align}\label{eq:example_s1_pvar}
s_1(x_1;\xi)=c_1^T \bar{Z}_2(x_1;\xi)=\underbrace{\left[{\footnotesize\mat{0\\0\\0\\1\\0\\0\\0\\2\\0\\0\\0\\1}}\right]^T}_{c_1^T}\underbrace{\left[{\footnotesize\mat{1\\x_1\\x_1^2\\\xi_1\\\xi_1 x_1\\\xi_1 x_1^2\\\xi_2\\\xi_2 x_1\\\xi_2 x_1^2\\\xi_3\\\xi_3 x_1\\\xi_3 x_1^2}}\right]}_{\bar{Z}_2(x_1;\xi)},
\end{align}
Here, the monomial bases $Z_2\in\R^{3}[x_1]$ and $\bar{Z}_2\in\R^{12}[x_1;\xi]$ are numerically represented by degree matrices $Z_{\text{M},2}\in\N^{3\times 1}$ and $\bar{Z}_{\text{M},2}\in\N^{12\times 4}$ respectively, defined as 	\vspace*{-0.3cm}
\begin{align}\label{eq:example_monomials}
Z_{\text{M},2}&:=\overbrace{\bmat{0\\1\\2}}^{x_1},	&	&\text{and} &
\bar{Z}_{\text{M},2}&:=\overbrace{\semismallbmat{0&0&0&0\\1&0&0&0\\2&0&0&0\\0&1&0&0\\1&1&0&0\\2&1&0&0\\0&0&1&0\\1&0&1&0\\2&0&1&0\\0&0&0&1\\1&0&0&1\\2&0&0&1}}^{[x_1,\xi_1,\xi_2,\xi_3]}
\end{align}

\section{Operations in the dpvar Representation}\label{sec:dpvar_operations}

We first show that, using the dpvar representation, standard operations on polynomial variables $s\in\R[x;\xi]$ may be performed at relatively low computational cost, by exploiting the affine contribution of the decision variables. In particular, we note that in the dpvar representation,
\begin{align*}
 s(x;\xi)&=Z_1(\xi)^T C Z_{d}(x)=\bmat{1\\\xi}^T C Z_{d}(x),
\end{align*}
so the vector of linear monomials $Z_1(\xi)$ always takes the same form. Therefore, there is no need to explicitly store or account for the degrees of the monomials in $Z_1(\xi)$, and the complexity of operations will be largely independent of the number of decision variables $\xi$.

By contrast, in the pvar representation,
\begin{align*}
 s(x;\xi)&=c^T \bar{Z}_d(x;\xi),
\end{align*}
the decision and independent variables are included in a single vector of monomials $\bar{Z}_d(x;\xi)$, taking the form
\begin{align}\label{eq:Zbar}
 \bar{Z}_d(x;\xi)=\bmat{1\\\xi}\otimes Z_{d}(x).
\end{align}
In this format, the decision variables and independent variables are represented using a single set of monomials. Implementing a data structure based on the pvar representation, therefore, the degrees of the decision variables $\xi$ have to be explicitly stored and processed when performing polynomial operations. As a result, the computational complexity of operating on the monomials will scale directly with the number of decision variables, even if the considered operation does not affect the decision variables (see Subsection~\ref{sec:subsec_miscellaneous_operations}). 

In the remainder of this section, we show how efficient addition, multiplication, and differentiation of polynomial variables may be performed using the dpvar representation. For each operation, the reduction in complexity using the dpvar representation is illustrated through a scalability test, comparing the time required to perform the operation using the \texttt{dpvar} data structure from SOSTOOLS 4.00, the \texttt{pvar} and \texttt{syms} structures from SOSTOOLS 3.04, as well as the YALMIP \texttt{sdpvar} structure. For the \texttt{syms} tests, the presented computation times include those necessary to convert the output to a (pvar) representation in terms of monomial degrees and coefficients, as needed for further processing in SOSTOOLS 3.04. All tests were performed on a computer with Intel Core i7-5960X CPU, and 128 GB of installed RAM.

\subsection{Addition}\label{sec:subsec_addition}

We first consider the operation of adding two (scalar) polynomial variables $s_1\in\R[x_1,\hdots,x_{p_1};\xi_1,\hdots,\xi_{q_1}]$ and $s_2\in\R[y_1,\hdots,y_{p_2};\eta_1,\hdots,\eta_{q_2}]$, written in the dpvar representation as
\begin{align*}
s_1(x;\xi)&=Z_1(\xi)^T C_1 Z_{d_1}(x),	&
s_2(y;\eta)&=Z_1(\eta)^T C_2 Z_{d_2}(y).
\end{align*}
In this format, it is clear that the sum $s_3=s_1+s_2$ of the polynomials may be expressed as
\begin{align*}
s_3(x,y;\xi,\eta)&=\bmat{Z_1(\xi)\\Z_1(\eta)}^T \bmat{C_1&0\\0&C_2} \bmat{Z_{d_1}(x)\\Z_{d_2}(y)}.
\end{align*}
The computational challenge, then, lies in defining the variables $z,\chi$, monomial basis $\hat{Z}_{d_3}\in\R^{n_3}[z]$, and coefficients $C_3$ to represent this result in the dpvar format,
\begin{align*}
s_3(z;\chi)&=Z_1(\chi)^T C_3 \hat{Z}_{d_3}(z)=\bmat{1\\\chi}^T C_3\hat{Z}_{d_3}(z).
\end{align*}
This may be achieved through the following steps:
\begin{enumerate}
	\item
	Combining the decision variables into a single vector $Z_1(\chi)$, where $\chi=\text{unique}(\xi;\eta)$.
	
	\item
	Combining the monomial bases $Z_{d_1}(x)$ and $Z_{d_2}(y)$ into a single vector $\hat{Z}_{d_3}(z)$, where $z=\text{unique}(x;y)$, and $d_3=\max\{d_1,d_2\}$.
	
	\item
	Rearranging and adding the elements of the coefficient matrix $\text{diag}(C_1,C_2)$ in accordance with the adjustments performed in the previous two steps.	
	
\end{enumerate}

Performing this conversion to the dpvar format, the greatest computational effort will generally be spent on the last two steps. Specifically, as shown in Appx.~\ref{apndx:merging_bases}, the complexity of merging degree matrices $Z_{\text{M},d_1}\in\N^{n_1\times p_1}$ and $Z_{\text{M},d_2}\in\N^{n_2\times p_2}$ is
\begin{align*}
	\mcl{O}\bl((n_1+n_2)\log(n_1+n_2)\br),
\end{align*}
where $n_i:=\frac{(p_i+d_i)!}{p_i!\ d_i!}$ denotes the number of monomials of degree at most $d_i$ in $p_i$ variables. For step 3, storing $C_1$ and $C_2$ as sparse matrices, the complexity of performing pre-established row and column permutations on $\text{diag}(C_1,C_2)$ will scale directly with the total number of nonzero coefficients as
\[
 \mcl{O}\bl(nnz(C_1)+nnz(C_2)\br),
\]
where the number of nonzero coefficients corresponds to the number of terms in each polynomial.
Notably, neither the complexity associated with step 2 nor that associated with step 3 depends directly on the number of decision variables, increasing only indirectly with the number of decision variables through the number of nonzero coefficients.

Now, compare this complexity to that of adding the same polynomials in the pvar representation,
\begin{align*}
s_1(x;\xi)&=c_1^T \bar{Z}_{d_1}(x;\xi),	&
s_2(y;\eta)&=c_2^T \bar{Z}_{d_2}(y;\eta),
\end{align*}
where $\bar{Z}_d$ is as in~\eqref{eq:Zbar}. Then
\begin{align*}
s_3(x,y;\xi,\eta)&=\bmat{c_1^T&c_2^T} \bmat{\bar{Z}_{d_1}(x;\xi)\\\bar{Z}_{d_2}(y;\eta)},
\end{align*}
once more requiring the monomial bases and coefficients to be combined.
In this case too, the complexity associated to combining the coefficients will scale as
\[
\mcl{O}\bl(nnz(c_1)+nnz(c_2)\br)=\mcl{O}\bl(nnz(C_1)+nnz(C_2)\br),
\]
requiring similar computational effort as when using the dpvar representation. However, since the number of monomials $\bar{n}_i$ in each vector $\bar{Z}_{d_i}$ now increases directly with the number of decision variables $q_i$ in each polynomial, 
\begin{align*}
	\bar{n}_i=(q_i+1)\cdot n_i = (q_i+1) \frac{(p_i+d_i)!}{p_i!\ d_i!}.
\end{align*}
the complexity of merging the bases will also increase with the number of decision variables,
\begin{align*}
	&\mcl{O}\bl((\bar{n}_1+\bar{n}_2)\log(\bar{n}_1+\bar{n}_2)\br)	\\
	&=\mcl{O}\bbl(\!\bl([q_1\!+\!1] n_1 + [q_2\!+\!1]  n_2\br)\log\bl([q_1\!+\!1] n_1 + [q_2\!+\!1] n_2\br)\!\bbr).
\end{align*}
For polynomials involving large numbers of decision variables $q_1$ and $q_2$, this complexity will be substantially worse than that of merging the bases in the dpvar representation.

\begin{figure*}[t!]
	\centering
	\begin{subfigure}[h]{0.49\textwidth}
		\centering
		\hspace*{-0cm}\includegraphics[width=1.0\textwidth]{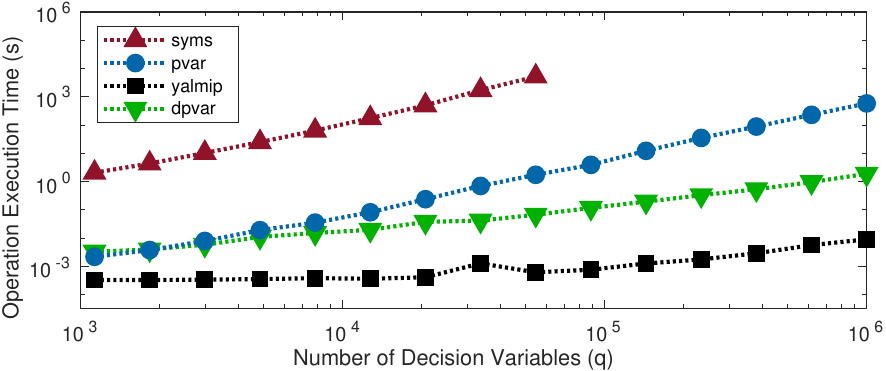}
		\caption{\footnotesize Computation time for addition $s_1(x;\xi)+s_2(y;\eta)$}
		\label{fig:addition}
	\end{subfigure}
	\begin{subfigure}[h]{0.49\textwidth}
		\centering
		\hspace*{-0cm}\includegraphics[width=1.0\textwidth]{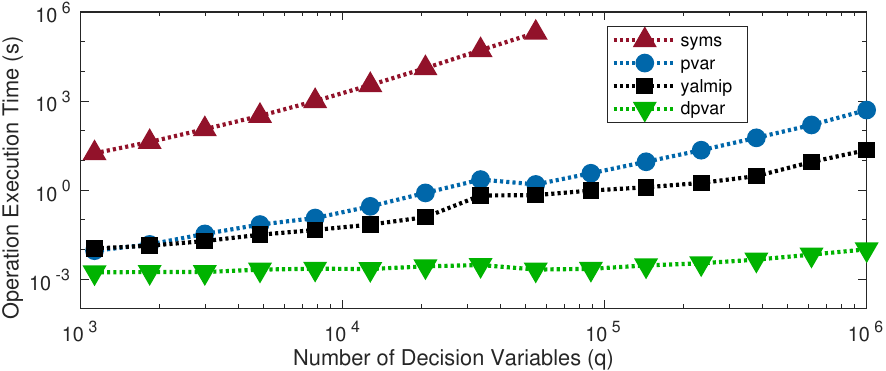}
		\caption{\footnotesize Computation time for multiplication $s_1(x;\xi)p_2(y)$}
		\label{fig:multiplication}
	\end{subfigure}
	
	\caption{\small Computation time for polynomial addition and multiplication using the \texttt{syms}, \texttt{pvar}, and \texttt{dpvar} data structures from respectively SOSTOOLS 3.04 and 4.00, and the \texttt{sdpvar} structure from YALMIP to represent the polynomials. The rate at which the computation time increases is relatively small using the \texttt{dpvar} structure compared to the alternatives, particularly for the multiplication operation. Only YALMIP achieves better performance for addition, by representing each monomial as a single index rather than as a set of degrees, requiring minimal computational effort to merge the bases of $s_1$ and $s_2$.
	}
	\label{fig:binary_operations}
	\vspace*{-0.4cm}
\end{figure*}

\paragraph*{\textbf{Example}} Consider the SOS variable $s_1(x_1;\xi):=\semismallbmat{1\\x_1}^T\semismallbmat{\xi_1&\xi_2\\\xi_2&\xi_3}\semismallbmat{1\\x_1}$. Defining $C_1\in\R^{4\times 3}$ as in Eqn.~\eqref{eq:example_s1_dpvar}, the sum $s_3(x_1;\xi)=s_1(x_1;\xi)+s_1(x_1;\xi)$ can then be represented in the dpvar format as
\begin{align*}
 s_3(x_1;\xi)=\bmat{Z_1(\xi)\\Z_1(\xi)}^T\bmat{C_1&0\\0&C_1}\bmat{Z_2(x_1)\\Z_2(x_1)}.
\end{align*}
Here, the computational cost of merging the decision variables is very small, and it is easy to recognize that the sum may be equivalently represented as
\begin{align*}
s_3(x_1;\xi)=Z_1(\xi)^T\bmat{C_1&C_1}\bmat{Z_2(x_1)\\Z_2(x_1)}.
\end{align*}
Similarly, it is computationally inexpensive to determine that the monomial vector $\hat{Z}_{2}(x_1)=\smallbmat{Z_2(x_1)\\Z_2(x_1)}$ pertains only a single independent variable $x_1$, and therefore, this vector may be numerically represented by the degree matrix
\begin{align*}
\hat{Z}_{\text{M},2}&=\overbrace{\bmat{Z_{\text{M},2}\\Z_{\text{M},2}}}^{[\ x_1 \ ]}\in\N^{6\times 1},
\end{align*}
where $Z_{\text{M},2}$ is as in Eqn.~\eqref{eq:example_monomials}.
Checking this matrix for unique monomials, only six rows have to be compared, and relatively little computational effort is necessary to establish a unique set of degrees, and to merge the columns of the coefficient matrix $\bmat{C_1&C_1}$ to find
\begin{align*}
s_3(x_1;\xi)=Z_1(\xi)^T[C_1+C_1]Z_2(x_1).
\end{align*}
Consider now computing the sum $s_3(x_1;\xi)=s_1(x_1;\xi)+s_1(x_1;\xi)$ using the pvar representation as
\begin{align*}
s_3(x_1;\xi)=\bmat{c_1^T&c_1^T}\bmat{\bar{Z}_2(x_1;\xi)\\\bar{Z}_2(x_1;\xi)}.
\end{align*}
where we define $c_1\in\R^{12}$ as in Eqn.~\eqref{eq:example_s1_pvar}. 
In this case, a unique set of variables $(x_1,\xi_1,\xi_2,\xi_3)$ can once again be established at relatively low computational cost, finding that the monomials $\check{Z}_2(x_1;\xi):=\smallbmat{\bar{Z}_{2}(x_1;\xi)\\\bar{Z}_{2}(x_1;\xi)}$ can be represented by the degree matrix \vspace*{-0.3cm}
\begin{align*}
\check{Z}_{\text{M},2}&=\overbrace{\bmat{\bar{Z}_{\text{M},2}\\\bar{Z}_{\text{M},2}}}^{[\ (x_1,\xi)\ ]}\in\N^{24\times 4},
\end{align*}
where $\bar{Z}_{\text{M},2}$ is as in Eqn.~\eqref{eq:example_monomials}.
However, the number of rows in this matrix is 4 times greater than that in the dpvar case, thus requiring a substantially greater computational effort to establish a unique set of degrees. This effect will be even worse for polynomial variables involving larger numbers of decision variables, offering a significant reduction in computation time using the \texttt{dpvar} data structure.
\medskip

The reduction in computation time offered by the dpvar representation is illustrated in Figure~\ref{fig:addition}, displaying the elapsed time for adding SOS variables $s_1(x_1,x_2;\xi_1,\hdots,\xi_q)$ and $s_2(y_1,y_2;\eta_{1},\hdots,\eta_{q})$ using the \texttt{dpvar}, \texttt{pvar}, \texttt{syms} and \texttt{sdpvar} (YALMIP) data structures, for increasing numbers of decision variables $q$. For each value of $q$, coefficients for $s_1$ and $s_2$ were randomly generated, and monomials $Z_{d}(x)$, $Z_{d}(y)$ of maximal degree $d=4$ were used. The decision variables were chosen such that $s_1$ and $s_2$ shared $\f{1}{2}q$ common variables, letting $\eta_j=\xi_{j+\frac{1}{2}q}$ for $j\in\{1,\hdots,\f{1}{2}q\}$.

\begin{figure*}[t!]
	\centering
	\hspace*{-0cm}\includegraphics[width=0.9\textwidth]{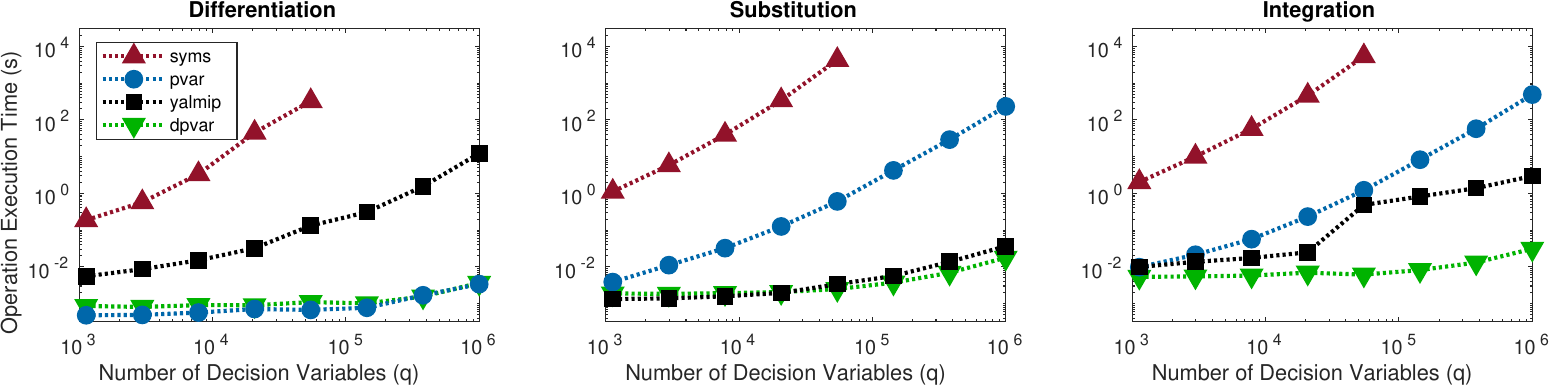}
	\caption{\small Computation time for differentiation, substitution, and integration of polynomial variables $s(x;\xi)$ using the \texttt{dpvar}, \texttt{pvar}, and \texttt{syms} data structures from SOSTOOLS 4.00, and the \texttt{sdpvar} data structure from YALMIP to represent $s$. Using the \texttt{dpvar} representation, the required time to perform each operation remains almost constant as the number of decision variables increases, offering substantial reductions in computation time for larger numbers of variables, compared to the alternative structures.}
	\label{fig:polynomial_operations}
	\vspace*{-0.4cm}
\end{figure*}

\subsection{Multiplication}\label{sec:subsec_multiplication}

We now consider the operation of polynomial multiplication, showing that this operation may also be performed more efficiently using the dpvar representation. For multiplication, since decision variables must always appear linearly in any SOS program, polynomial variables $s\in\R[x;\xi]$ may only be multiplied by known polynomial functions $p\in\R[y]$. In dpvar format, these may be expressed as
\begin{align*}
s_1(x;\xi)&=Z_1(\xi)^T C Z_{d_1}(x),	&
p_2(y)&= b^T Z_{d_2}(y),
\end{align*}
so that the product becomes
\begin{align*}
s_1(x;\xi)p_2(y)=Z_1(\xi)^T \bl(b^T\otimes C\br)\bl(Z_{d_2}(y)\otimes Z_{d_1}(x)\br).
\end{align*}
Performing this operation in MATLAB, the coefficients $b,C$ and monomial degrees $Z_{\text{M},d_1}(x),Z_{\text{M},d_2}(y)$ may be stored as sparse matrices. Then, performing the Kronecker product $b^T\otimes C$ will require multiplying at most $nnz(C)\cdot nnz(b)$ elements, invoking a worst-case complexity of
\[
\mcl{O}\bl(nnz(C)nnz(b)\br).
\]
To compute the product $Z_{d_2}(y)\otimes Z_{d_1}(x)$, the nonzero degrees of all the variables in each monomial in $Z_{d_2}$ must be added to the degrees of the same variables in each of the monomials in $Z_{d_1}$. In the worst-case scenario (e.g. $x=y$ and $Z_{d_1}=Z_{d_2}$), this will require adding all nonzero degrees in $Z_{\text{M},d_2}$ to all nonzero degrees in $Z_{\text{M},d_1}$. The complexity of this operation scales as
\[
 \mcl{O}\bl(nnz(Z_{\text{M},d_1})nnz(Z_{\text{M},d_2})\br).
\]
Consider now computing the same product based on the pvar representation,
\begin{align*}
s_1(x;\xi)&=c^T \bar{Z}_{d_1}(x;\xi),	&
p_2(y)&= b^T Z_{d_2}(y),
\end{align*}
so that
\begin{align*}
s_1(x;\xi)p_2(y)=(b^T\otimes c^T)(Z_{d_2}(y)\otimes \bar{Z}_{d_1}(x;\xi)).
\end{align*}
As was the case in the dpvar representation, the cost of computing the new coefficients will be
\[
 \mcl{O}\bl(nnz(c)nnz(b)\br)=\mcl{O}\bl(nnz(C)nnz(b)\br),
\]
scaling with the product of the number of terms in the two polynomials. However, in the pvar representation, the number of nonzero degrees in $\bar{Z}_{\text{M},d_1}$ increases linearly with the number of decision variables $q$ in $s_1$, so that the complexity of multiplying the bases will be
\[
 \mcl{O}\bl(nnz(\bar{Z}_{\text{M},d_1})nnz(Z_{\text{M},d_2})\br)
 =\mcl{O}\bl(q\cdot nnz(Z_{\text{M},d_1})nnz(Z_{\text{M},d_2})\br).
\]
This dependence on the number of decision variables is not present when implementing the dpvar representation, resulting in a substantial difference in computational complexity for large values of $q$. 

\paragraph*{\textbf{Example}} Consider the polynomial function $p_1(x_1)=1-2x_1^2$ and the SOS variable $s_1(x_1;\xi):=\semismallbmat{1\\x_1}^T\semismallbmat{\xi_1&\xi_2\\\xi_2&\xi_3}\semismallbmat{1\\x_1}$. Defining $b_1\in\R^{3}$ as in Eqn.~\eqref{eq:example_p1_pvar} and $C_1\in\R^{4\times 3}$ as in Eqn.~\eqref{eq:example_s1_dpvar}, the product $s_3(x_1;\xi)=s_1(x_1;\xi)p_1(x_1)$ can then be represented in the dpvar format as
\begin{align*}
 s_3(x_1;\xi)&=Z_1(\xi)^T\bl(b_1^T\otimes C_1\br)\bl(Z_{2}(x_1)\otimes Z_{2}(x_1)\br).
\end{align*}
Similarly, defining $c_1\in\R^{12}$ as in Eqn.~\eqref{eq:example_s1_pvar}, the product $s_3(x_1;\xi)=s_1(x_1;\xi)p_1(x_1)$ can also be represented in the pvar format as
\begin{align*}
s_3(x_1;\xi)&=\bl(b_1^T\otimes c_1^T\br)\bl(Z_{2}(x_1)\otimes \bar{Z}_{2}(x_1;\xi)\br).
\end{align*}
Here, the monomial vectors $Z_{2}\in\R^{3}[x_1]$ and $\bar{Z}_{2}\in\R^{12}[x_1;\xi]$ can be represented by respectively the degree matrix $Z_{\text{M},2}\in\N^{3\times 1}$ and $Z_{\text{M},2}\in\N^{12\times 4}$ as in Eqn.~\eqref{eq:example_monomials}. However, where the former degree matrix contains only 2 nonzero elements, the latter matrix contains $17$ nonzero elements. As such, the cost of computing the degree matrix associated to the Kronecker product $Z_{2}(x_1)\otimes \bar{Z}_{2}(x_1;\xi)$ will also be more than 8 times as great as that of computing the degrees for $Z_{2}(x_1)\otimes Z_{2}(x_1)$.
\medskip

The reduction in complexity offered by the dpvar representation can also be observed in Figure~\ref{fig:multiplication}, displaying the elapsed time for multiplying a randomly generated variable $s_1(x_1,x_2;\xi_1,\hdots,\xi_q)$ (see Subsection~\ref{sec:subsec_addition}) and polynomial $p_2(y_1,y_2)$ using the different data structures.

\subsection{Differentiation, Substitution, and Integration}\label{sec:subsec_miscellaneous_operations}

Finally, we consider the operations of differentiation, substitution and integration. For an arbitrary polynomial $s\in\R[x;\xi]$ in the dpvar representation,
\begin{align*}
s(x;\xi)&=Z_1(\xi)^T C Z_{d}(x),
\end{align*}
these operations will involve only adjusting the monomial vector $Z_{d}$, and associated columns in the coefficient matrix $C$. For example, let $z_{ij}=[Z_{d}]_{ij}$ denote the element in row $i$ and column $j$ of the degree matrix $Z_{\text{M},d}\in\N^{n\times p}$, and let $C_i$ denote the $i$th column of the coefficient matrix $C\in\R^{(q+1)\times n}$. Then, differentiation with respect to $x_j$ may be performed by multiplying all elements in each column $C_i$ for $i=1,\hdots,n$ with $z_{ij}$, and subtracting a value of $1$ from all nonzero degrees in column $j$ of $Z_{\text{M},d}\in\N^{n\times p}$.
The complexity of this operation depends only indirectly on the number of decision variables, as each decision variable adds a row to the coefficient matrix $C\in\R^{(q+1)\times n}$.

By contrast, performing the same operations using the pvar representation, 
\begin{align*}
s(x;\xi)&=b^T \bar{Z}_{d}(x;\xi),
\end{align*}
the decision variables are included in the monomial basis $\bar{Z}_d$. Therefore, the complexity of finding and adjusting the appropriate degrees of the monomials to account for e.g. differentiation with respect to a variable $x_j$, will directly increase with the number of decision variables, despite the fact that the decision variables themselves are invariant under these operations. In this sense, unnecessary computational overhead is introduced when performing differentiation, substitution and integration in the pvar representation, which is avoided implementing the dpvar representation. 

\paragraph*{\textbf{Example}} Consider the SOS variable $s_1(x_1;\xi):=\semismallbmat{1\\x_1}^T\semismallbmat{\xi_1&\xi_2\\\xi_2&\xi_3}\semismallbmat{1\\x_1}$, represented in the dpvar representation as
\begin{align*}
 s_1(x_1;\xi)=Z_1(\xi)^T C_1 Z_2(x_1) =\semismallbmat{1\\\xi_1\\\xi_2\\\xi_3}^T\semismallbmat{0&0&0\\1&0&0\\0&2&0\\0&0&1}\semismallbmat{1\\ x_1\\ x_1^2}
\end{align*}
Then the derivative of this variable with respect to $x_1$ can be easily obtained by multiplying each column in $C_1$ with their associated degree in $Z_{\text{M},2}\in\N^{3\times 1}$, and reducing all nonzero degrees with a value of 1:
\begin{align*}
 \f{\partial}{\partial x_1}s(x_1)
 =\semismallbmat{1\\\xi_1\\\xi_2\\\xi_3}^T\semismallbmat{0&0&0\\0&0&0\\0&2&0\\0&0&2}\semismallbmat{1\\1\\ x_1}
 =\semismallbmat{1\\\xi_1\\\xi_2\\\xi_3}^T\semismallbmat{0&0\\0&0\\2&0\\0&2}\semismallbmat{1\\ x_1}.
\end{align*}
Numerically, this requires only multiplying two nonzero degrees with two nonzero coefficients, and then subtracting a value of 1 from these two nonzero degrees.
By contrast, in the pvar representation, \vspace*{-0.3cm}
\begin{align*}
s_1(x_1;\xi)=c_1^T \bar{Z}_2(x_1;\xi)=\left[{\footnotesize\mat{0\\0\\0\\1\\0\\0\\0\\2\\0\\0\\0\\1}}\right]^T\left[{\footnotesize\mat{1\\x_1\\x_1^2\\\xi_1\\\xi_1 x_1\\\xi_1 x_1^2\\\xi_2\\\xi_2 x_1\\\xi_2 x_1^2\\\xi_3\\\xi_3 x_1\\\xi_3 x_1^2}}\right],
\end{align*}
the degree matrix $\bar{Z}_{\text{M},2}\in\N^{24\times 4}$ has eight nonzero elements in the column associated to the variable $x_1$. Although the computational cost of subtracting a value of 1 from each of these degrees will not be substantial in this case, for examples involving larger numbers of decision variables, this may amount to a nontrivial reduction in computational complexity using the dpvar representation.
\medskip

The reduced computation time allowed by the dpvar representation for larger-scale tests is illustrated in Figure~\ref{fig:polynomial_operations}, presenting the elapsed time for differentiation, substitution and integration of a randomly generated polynomial $s_1(x_1,x_2;\xi_1,\hdots,\xi_q)$ with respect to the variable $x_2$, using the different SOSTOOLS and YALMIP data structures, and for increasing numbers of decision variables $q$.

\vspace*{-0.1cm}

\section{Storage and Manipulation of dpvars}\label{sec:dpvar_storage}

Having analyzed the complexity of standard operations in the dpvar representation, in this section, we show how this representation also allows the memory burden and general computational overhead that comes with parsing an SOS program to be reduced. In particular, implementing the dpvar representation in MATLAB, we define a polynomial variable $S\in\R^{m_1\times m_2}[x;\xi]$ using the \texttt{dpvar} data structure, storing
\begin{itemize}
	\item The independent variables $x_1,\hdots,x_p$.	
	\item The decision variables $\xi_1,\hdots,\xi_q$.
	\item The monomial degrees $Z_{\text{M},d}\in\N^{n\times p}$.
	\item The coefficient matrix $C\in\R^{m_1(q+1)\times m_2 n}$.
\end{itemize}
Decomposing the polynomial in this manner, the greatest storage cost will be that associated to the monomial degrees $Z_{\text{M},d}$ and coefficient matrix $C$. However, storing both of these fields as sparse matrices in MATLAB, the memory overhead will be minimal, as we show in Subsection~\ref{sec:subsec_linearity_in_memory}.
In addition, exploiting the structure of \texttt{dpvar} objects, matrix operations such as concatenation can be performed with relatively low computational overhead, as detailed in Subsection~\ref{sec:subsec_dpvar_manipulation}.
\vspace*{-0.2cm}

\subsection{Memory Complexity of Storing dpvar Objects}\label{sec:subsec_linearity_in_memory}

Exploiting linearity of the decision variables in its structure, the dpvar representation allows polynomial variables to be stored in programming languages with sparsity structures using minimal memory with respect to the number of decision variables.
Specifically, consider storing a matrix-valued polynomial variable $S\in\R^{m_1\times m_2}[x_1,\hdots,x_p;\xi_1,\hdots,\xi_q]$, expressed in the dpvar representation as
\begin{align}\label{eq:dpvar_format}
S(x;\xi)&=\bl(I_{m_1}\otimes Z_1(\xi)\br)^T C \bl(I_{m_2}\otimes Z_{d}(x)\br).
\end{align}
As mentioned, the greatest memory burden in representing this variable in MATLAB will be that associated to storing the coefficient matrix $C\in\R^{m_1 (q+1)\times m_2 n_1}$, and the monomial degrees $Z_{\text{M},d}\in\N^{n_1\times p}$. Storing both objects as sparse matrices, only the nonzero coefficients and degrees are retained, so that the required memory scales as
\begin{align*}
\mcl{O}\bl(nnz(C)+nnz(Z_{\text{M},d})\br).
\end{align*}
This cost does not depend directly on the number of decision variables.

Consider now storing the same variable in the pvar format,
\begin{align}\label{eq:pvar_format}
 S(x;\xi):=B^T \bl(I_{m_2}\otimes \bar{Z}_d(x;\xi)\br),
\end{align}
where $B\in\R^{m_1\times m_2n_2}$ and $\bar{Z}_d=\bmat{1\\\xi}\otimes Z_{d}(x)\in\R^{n_2}[x;\xi]$. Using this representation, the storage cost will also mostly be determined by the number of nonzero coefficients and degrees. Since the number of nonzero coefficients is independent of the representation, the cost of storing these coefficients will be roughly the same using the dpvar and pvar structures, scaling with $nnz(C)=nnz(B)$. However, when considering $q$ decision variables, each monomial appearing in the vector $Z_{d}(x)$ will appear $q+1$ times in the vector $\bar{Z}_d(x;\xi)$. Therefore, each nonzero degree in $Z_{\text{M},d}\in\N^{n_1\times p}$ will also appear $q+1$ times in $\bar{Z}_{\text{M},d}\in\N^{(q+1)n_1\times (q+p)}$. Moreover, for each of the $n_1$ monomials included in $Z_{d}(x)$, the nonzero degrees of the decision variables will also need to be stored, amounting to a total number of $nnz(\bar{Z}_{\text{M},d}) = (q+1)nnz(Z_{\text{M},d})+q n_1$ nonzero degrees,
\begin{align*}
 nnz(\bar{Z}_{\text{M},d}) = (q+1)nnz(Z_{\text{M},d})+q n_1
\end{align*}
The cost of storing the coefficients and monomials in the pvar representation thus scales with
\begin{align*}
\mcl{O}\bl(nnz(C)+(q+1)nnz(Z_{\text{M},d})+qn_1\br).
\end{align*}
Implementing the pvar representation, the required memory of storing the monomials increases directly with the number of decision variables. 
For large numbers of decision variables $q$, this amounts to a substantial storage cost that may be avoided using the \texttt{dpvar} structure.

\paragraph*{\textbf{Example}} Numerically representing the SOS variable $s_1(x_1;\xi):=\semismallbmat{1\\x_1}^T\semismallbmat{\xi_1&\xi_2\\\xi_2&\xi_3}\semismallbmat{1\\x_1}$ in the dpvar format (Eqn.~\eqref{eq:example_s1_dpvar}), only 2 nonzero degrees have to be stored. By contrast, representing this variable in the pvar format (Eqn.~\eqref{eq:example_s1_pvar}), 17 nonzero degrees have to be stored. Including the 3 nonzero coefficients in each representation, the total number of nonzero elements that need to be stored to represent $s_1$ is 4 times smaller using the \texttt{dpvar} structure than using the \texttt{pvar} structure (see also Section~\ref{sec:MATLAB_sparsity}).

\vspace*{-0.2cm}

\subsection{Matrix Operations on \texttt{dpvar} Objects}\label{sec:subsec_dpvar_manipulation}

In many SOS programs, the polynomial decision variables appear as matrix-valued objects. Therefore, in addition to the standard polynomial operations discussed in Section~\ref{sec:dpvar_operations}, matrix operations such as concatenation must also be efficiently implemented in any SOS parser. Using the dpvar representation, this can be achieved by exploiting the block structure of the coefficient matrix. In particular, for a variable $S\in\R^{m_1\times m_2}[x;\xi]$, the coefficient matrix $C\in\R^{m_1(q+1)\times m_2 n}$ is comprised of $m_1\times m_2$ blocks $C_{ij}\in\R^{(q+1)\times n}$, each corresponding to a single element of the matrix-valued variable. This allows for efficient assignment and modification of individual elements of the polynomial variable. In addition, for two matrix-valued polynomial variables $S_1,S_2\in\R^{m_1\times m_2}[x;\xi]$, defined in terms of the same monomial basis $Z_{d}$ as 
\begin{align*}
S_i(x;\xi)&=\bl(I_{m_1}\otimes Z_1(\xi)\br)^T C_i \bl(I_{m_2}\otimes Z_{d}(x)\br),
\end{align*}
concatenation of $S_1$ and $S_2$ merely requires concatenating the coefficient matrices $C_1$ and $C_2$. For example, vertical concatenation of $S_1,S_2$ may be represented as
\begin{align*}
&\bmat{S_1(x;\xi)\\S_2(x;\xi)}	
=
\bl(I_{2m_1}\otimes Z_1(\xi)\br)^T\bmat{C_1\\ C_2}\bl(I_{m_2}\otimes Z_{d}(x)\br),
\end{align*}
requiring almost no computational effort.
Of course, if $S_1$ and $S_2$ are defined in terms of different monomial bases, these bases would have to be merged first, for which we refer to the discussion in Subsection~\ref{sec:subsec_addition}.

\section{Exploiting Sparsity in Storage and Operation}\label{sec:MATLAB_sparsity}

Having presented the benefits of using the dpvar representation in parsing SOS programs, we finally show how the \texttt{dpvar} data structure exploits the MATLAB built-in sparsity structure to minimize memory and computational overhead in numerically representing polynomial variables. In particular, in Subsection~\ref{sec:subsec_CSC_format}, we outline how sparse matrices are implemented in MATLAB and analyze how this format affects memory and computational complexity. In Subsection~\ref{sec:subsec_sparsity_dpvar}, we subsequently show how the \texttt{dpvar} data structure exploits this format in storing the coefficient matrix and monomial degrees, to optimize performance.

\subsection{The Compressed Sparse Column Format}\label{sec:subsec_CSC_format}

In MATLAB, the built-in sparse storage structure is optimized for storing and operating on matrices with relatively few columns. In particular, sparse matrices are implemented using a Compressed Sparse Column (CSC) format~\cite{gilbert1992SparseMatlab}, representing a matrix $A\in\R^{m\times n}$ with $nnz(A)$ nonzero elements through three arrays:
\begin{enumerate}
	\item An array $\texttt{a}\in\R^{nnz(A)}$ of nonzero elements.
	\item An array $\texttt{r}\in\R^{nnz(A)}$ of row indices.
	\item An array $\texttt{cp}\in\R^{n+1}$ of column pointers.
\end{enumerate}
In the first of these arrays, $\texttt{a}\in\R^{nnz(A)}$, all nonzero elements of the matrix are collected in \textit{column-major} order. That is, letting $\{a_1,\hdots,a_n\}$ denote the columns of the matrix $A$, and letting $\{\bar{a}_1,\hdots,\bar{a}_n\}$ denote the nonzero elements from these columns, the first array $\texttt{a}$ may be constructed as:
\[
 \texttt{a}=\bmat{\bar{a}_1^T,&\hdots,&\bar{a}_n^T}^T\in\R^{nnz(A)}.
\]
Corresponding row numbers for these nonzero elements are then stored in the array \texttt{r}, so that the $k$th nonzero element $\texttt{a}(k)$ appears in row $\texttt{r}(k)$ of the matrix $A$.
Finally, for each of the columns $j=1,\hdots,n$ of the matrix, a column pointer is stored in the array \texttt{cp}.  
Letting $\ell_j=nnz(a_j)$,
this column pointer is defined as
\begin{align*}
\texttt{cp}=\bmat{1,\ 1+\ell_1,\ \hdots,\ 1+\sum_{j=1}^{n-1}\ell_j,\  \sum_{j=1}^{n}\ell_j}\in\R^{n+1},
\end{align*}
so that $\texttt{a}\bl(\texttt{cp}(j)\br)$ provides the first nonzero element of column $j\in\{1,\hdots,n\}$ of $A\in\R^{m\times n}$. 

Using this data structure to store (sparse) matrices, the required memory will be minimal for matrices with few columns. In particular, although the cost of storing $\texttt{a}\in\R^{nnz(A)}$ and $\texttt{r}\in\R^{nnz(A)}$ depends only on the number of nonzero elements $nnz(A)$, the memory necessary to store the array $\texttt{cp}\in\R^{n+1}$ is determined by the number of columns $n$ of the matrix. Therefore, the memory burden for storing sparse matrices increases with the number of columns in this matrix, even if these columns do not contain any nonzero elements.

In addition, using the CSC storage format, the complexity of operations involving full or partial columns of the matrix will generally be smaller than those involving full or partial rows of the matrix. Indeed, for any column $j\in\{1,\hdots,n\}$ of $A$, the nonzero elements appearing in this column are known to be stored at positions $k\in\{\texttt{cp}(j),\texttt{cp}(j)+1,\hdots,\texttt{cp}(j+1)-1\}$ within the array $\texttt{a}$, requiring minimal effort to access these elements. On the other hand, in order to access elements of a particular row $i\in\{1,\hdots,m\}$ of the matrix, all indices $k\in\{1,\hdots,nnz(A)\}$ with associated row index $\texttt{r}(k)=i$ have to be found, potentially requiring the full array $\texttt{r}$ to be analyzed.
This introduces additional computational overhead when operating on full or partial rows of the matrix, generally making ``row-based'' operations more computationally demanding than ``column-based'' equivalents.

\subsection{Sparsity in the dpvar Structure}\label{sec:subsec_sparsity_dpvar}

We now show how, using the \texttt{dpvar} data structure, the CSC storage format may be exploited to minimize the storage and operational cost of representing and manipulating polynomial variables. To illustrate, consider storing a variable
\begin{align*}
 s(x;\xi)=Z_1(\xi)C Z_{d}(x)\ \in \R[x_1,\hdots x_p;\xi_1,\hdots\xi_q].
\end{align*}
Storing the coefficient matrix $C\in\R^{(q+1)\times n}$ and monomial degrees $Z_{\text{M},d}\in\N^{n\times p}$ using the CSC structure, the required memory will be relatively small. In particular, since $p$ variables allow $n=\frac{(p+d)!}{p!d!}$ monomials of degree at most $d$, the number of rows in the monomial degree matrix $Z_{\text{M},d}\in\N^{n\times p}$ will in general vastly exceed the number of columns. In addition, in SOS programs, a monomial $[Z_{d}]_k$ is often paired with multiple decision variables $\xi_j$. As a consequence, the number of decision variables tends to exceed the number of monomials, and thus the number of rows in the coefficient matrix $C\in\R^{(q+1)\times n}$ also tends to be at least as large as the number of columns. 
Since the memory cost of storing a matrix in the CSC format increases with the number of columns, the fact that both the coefficient matrix and monomial degree table contain relatively few columns allows polynomial variables to be efficiently stored using the \texttt{dpvar} data structure.

Similarly, the complexity of performing operations on variables in the \texttt{dpvar} structure may be minimized using the sparse storage structure. In particular, as discussed in Subsection~\ref{sec:subsec_addition}, a significant part of the computational complexity in performing operations such as addition comes from having to merge the rows of the monomial degree matrix $Z_{\text{M},d}\in\R^{n\times p}$, and associated columns of the coefficient matrix $C\in\R^{(q+1)\times n}$. Here, although the CSC storage format is poorly-suited for comparing the large amounts of rows in the monomial matrix, the small number of columns in $Z_{\text{M},d}$ ensures the complexity of this process remains relatively small. Moreover, the column-major storage structure allows the columns of the coefficient matrix to be permuted with relatively high efficiency, invoking a complexity that does not depend directly on the number of rows $(q+1)$ of $C$. 
Thus, exploiting the MATLAB sparse storage structure, the \texttt{dpvar} data structure allows the computational cost of operations like addition to be minimized with respect to the number of decision variables $q$.

\paragraph*{\textbf{Example}} Consider the SOS variable $s_1(x_1;\xi):=\semismallbmat{1\\x_1}^T\semismallbmat{\xi_1&\xi_2\\\xi_2&\xi_3}\semismallbmat{1\\x_1}$, which can be represented in the dpvar format as
\vspace*{-0.3cm}
\begin{align*}
s_1(x_1;\xi) &=Z_1(\xi)^T C_1 Z_2(x_1)
={\footnotesize\bmat{1\\\xi_1\\\xi_2\\\xi_3}}^T{\footnotesize\bmat{0&0&0\\1&0&0\\0&2&0\\0&0&1}}{\footnotesize\bmat{1\\ x_1\\x_1^2}},
\end{align*}
and in the pvar format as \vspace*{-0.3cm}
\begin{align*}
s_1(x_1;\xi)=c_1^T \bar{Z}_2(x_1;\xi)=\left[{\footnotesize\mat{0\\0\\0\\1\\0\\0\\0\\2\\0\\0\\0\\1}}\right]^T\left[{\footnotesize\mat{1\\x_1\\x_1^2\\\xi_1\\\xi_1 x_1\\\xi_1 x_1^2\\\xi_2\\\xi_2 x_1\\\xi_2 x_1^2\\\xi_3\\\xi_3 x_1\\\xi_3 x_1^2}}\right],
\end{align*}
where, the monomial bases $Z_2\in\R^{3}[x_1]$ and $\bar{Z}_2\in\R^{12}[x_1;\xi]$ are numerically represented by matrices	\vspace*{-0.3cm}
\begin{align*}
Z_{\text{M},2}&=\overbrace{\bmat{0\\1\\2}}^{x_1},	&	&\text{and} &
\bar{Z}_{\text{M},2}&=\overbrace{\semismallbmat{0&0&0&0\\1&0&0&0\\2&0&0&0\\0&1&0&0\\1&1&0&0\\2&1&0&0\\0&0&1&0\\1&0&1&0\\2&0&1&0\\0&0&0&1\\1&0&0&1\\2&0&0&1}}^{[x_1,\xi_1,\xi_2,\xi_3]}.
\end{align*}
Then, in the dpvar format, the coefficients $C_1$ can be stored in the CSC format as
\begin{align*}
\texttt{a}_{C_1}&=\bmat{1\\2\\1},	&
\texttt{r}_{C_1}&=\bmat{2\\3\\4},	&
\texttt{cp}_{C_1}&=\bmat{1\\2\\3\\3},
\end{align*}
where \texttt{a} denotes the array of nonzero elements, \texttt{r} the array of row numbers, and \texttt{cp} the array of column pointers. Similarly, the degree matrix $Z_{\text{M,2}}$ can be stored in the CSC format as
\begin{align*}
\texttt{a}_{Z_2}&=\bmat{1\\2},	&
\texttt{r}_{Z_2}&=\bmat{2\\3},	&
\texttt{cp}_{Z_2}&=\bmat{1\\2},
\end{align*}
requiring a total of 16 values to be stored in order to represent the coefficients and degrees using the \texttt{dpvar} structure. On the other hand, using the \texttt{pvar} structure, the coefficients $c_1$ are stored in CSC format as
\begin{align*}
\texttt{a}_{c_1}&=\bmat{1\\2\\1},	&
\texttt{r}_{c_1}&=\bmat{4\\8\\12},	&
\texttt{cp}_{C_1}&=\bmat{1\\3},
\end{align*}
and the degrees $\bar{Z}_{\text{M},2}$ are stored as
\begin{align*}
\texttt{a}_{\bar{Z}_2}&=\semismallbmat{1\\2\\1\\2\\1\\2\\1\\2\\1\\1\\1\\1\\1\\1\\1\\1\\1},	&
\texttt{r}_{\bar{Z}_2}&=\semismallbmat{2\\3\\5\\6\\8\\9\\11\\12\\4\\5\\6\\7\\8\\9\\10\\11\\12},	&
\texttt{cp}_{\bar{Z}_2}&=\semismallbmat{1\\9\\12\\15\\17}.	\\
\end{align*}
Although the \texttt{pvar} structure allows the coefficients to be stored slightly more efficiently, the memory required to store the degrees will be substantially larger, amounting to a total of 47 values to be stored to represent both the degrees and coefficients. This is almost 3 times as many values as using the \texttt{dpvar} structure, exemplifying the significant reduction in memory requirements that the \texttt{dpvar} structure allows.

\section{Incorporation into SOSTOOLS}\label{sec:SOSTOOLS_implementation}

Having demonstrated the advantages of using the \texttt{dpvar} data structure for parsing polynomial variables, we now consider the incorporation of this structure in SOSTOOLS. Specifically, for SOSTOOLS version 4.00~\cite{papachristodoulou2021SOSTOOLS}, we have modified all functions to use the \texttt{dpvar} data structure for definition and manipulation of (polynomial) decision variables. To illustrate the enhanced performance this offers, in this section, we consider several polynomial optimization problems that are commonly solved with SOSTOOLS. For each problem, we compare the time required for parsing the problem using SOSTOOLS 3.04, SOSTOOLS 4.00, and using the \textit{batch} parser YALMIP~\cite{lofberg2004YALMIP}. To solve the resulting SDP, in each case, SEDUMI~\cite{sturm1999SEDUMI} was used. More details on the exact implementation of each problem in SOSTOOLS may be found in Appx.~\ref{apndx:optimization_tests}.

\begin{figure*}[t!]
	\centering
	\begin{subfigure}[h]{0.32\textwidth}
		\centering
		\hspace*{-0cm}\includegraphics[width=1.0\textwidth]{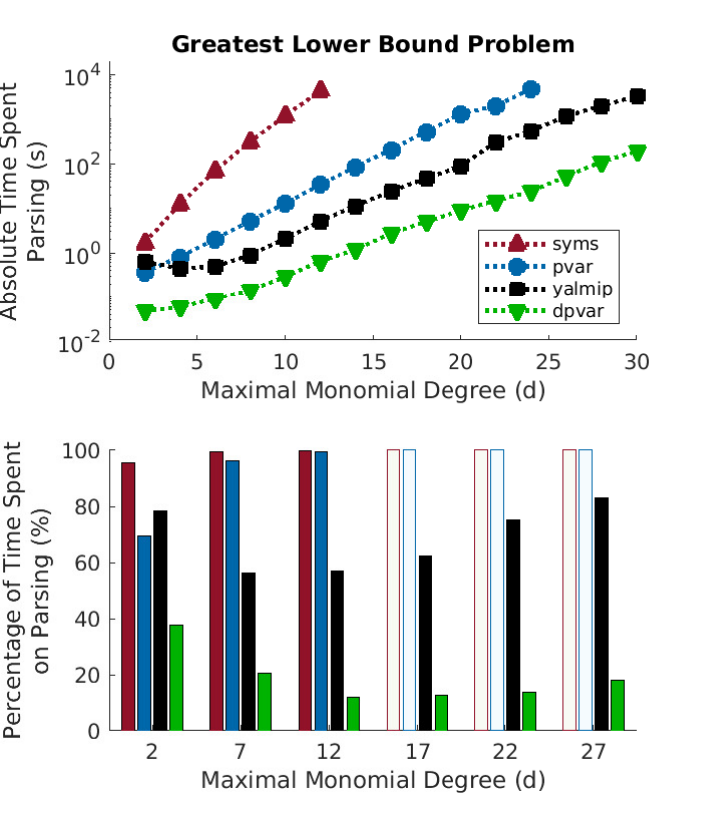}
		\vspace*{-0.6cm}
		\caption{\footnotesize \vspace*{-0.075cm} Greatest lower bound test, Subsection~\ref{sec:subsec_Greatest_Lower_Bound}}
		\label{fig:polynomial_optimization}
	\end{subfigure}
	\begin{subfigure}[h]{0.32\textwidth}
		\centering
		\hspace*{-0cm}\includegraphics[width=1.0\textwidth]{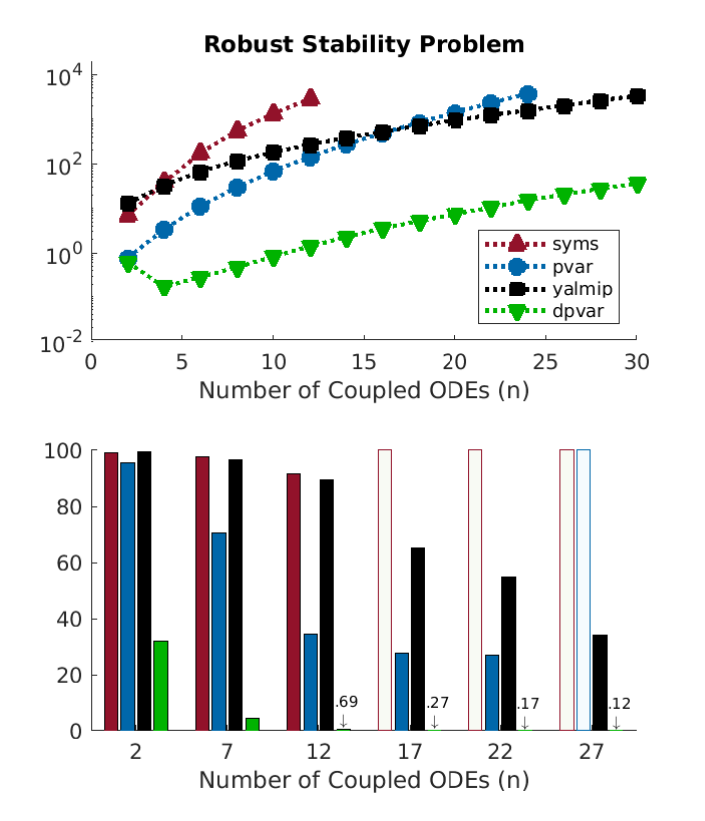}
		\vspace*{-0.6cm}
		\caption{\footnotesize \vspace*{-0.075cm} Robust stability test, Subsection~\ref{sec:subsec_Robust_Stability}}
		\label{fig:robust_stability}
	\end{subfigure}
	\begin{subfigure}[h]{0.32\textwidth}
		\centering
		\hspace*{-0cm}\includegraphics[width=1.0\textwidth]{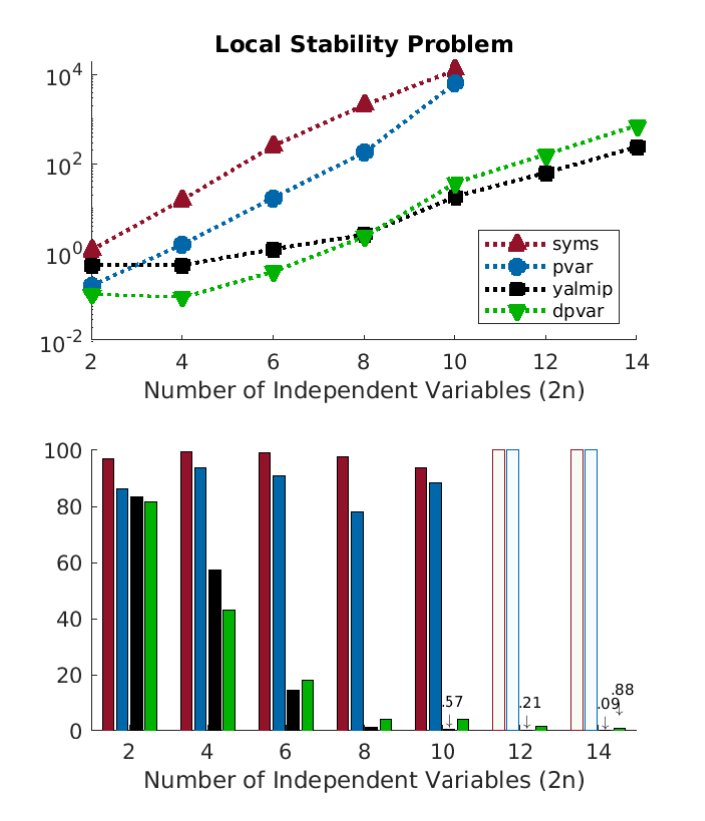}
		\vspace*{-0.6cm}
		\caption{\footnotesize \vspace*{-0.075cm} Local stability test, Subsection~\ref{sec:subsec_Local_Stability}}
		\label{fig:local_stability}
	\end{subfigure}
	
	\caption{\small Elapsed time parsing the polynomial optimization problems from Section~\ref{sec:SOSTOOLS_implementation}, using SOSTOOLS 3.04 with the \texttt{syms} and \texttt{pvar} data structures, SOSTOOLS 4.00 with the \texttt{dpvar} data structure, and using the batch parser YALMIP. Tests in each case were discontinued when the parsing time exceeded 3600 seconds, or the solver ran out of memory. The percentage of time spent parsing each problem was computed by dividing the absolute time spent parsing the SOS program by the sum of the time spent parsing the SOSP and solving the resulting SDP, for each implementation. The results show that, using the \texttt{dpvar} data structure, SOSTOOLS 4.00 is able to parse common SOS problems with an efficiency comparable to, or even greater than that using the batch parser YALMIP. 
	}
	\label{fig:SOSTOOLS_implementation_scalability_tests}
	\vspace*{-0.6cm}
\end{figure*}

\subsection{Greatest Lower Bound}\label{sec:subsec_Greatest_Lower_Bound}

As a first problem, we seek the greatest lower bound (GLB) $\gamma$ on some function $f$,
\begin{align*}
\max_{\gamma}\quad \gamma,&	&
\text{s.t.}\quad \gamma&\leq f(x) \quad
\forall x_1,x_2\in[-12,12],
\end{align*}
where $f(x)=x_1^4+x_2^4-2x_2x_1^3-3x_2^2x_1^2+150(x_1^2+x_2^2)$. To enforce the constraints $x_1,x_2\in[-12,12]$, we require
\begin{align*}
g_1(x)&=12^2-x_1^2\geq 0,	&
g_2(x)&=12^2-x_2^2\geq 0,\\
g_3(x)&=2\cdot 12^2-(x_1^2+x_2^2)\geq 0.
\end{align*}
Invoking Putinar's Positivstellensatz~\cite{putinar1993Psatz} (Psatz), we enforce a single SOS constraint
\begin{align}\label{eq:GLB_SOS}
(f(x)-\gamma) - s_1g_1(x) - s_2g_2(x) - s_3g_3(x)\in\Sigma_s,
\end{align}
with SOS variables $s_1,s_2,s_3\in\Sigma_s$.

In parsing the GLB program, the maximal degree of monomials $d$ appearing in the variables $s_i=Z_{d}(x)^T P Z_{d}(x)$ may be increased, allowing for more accurate results at the expense of a higher computational complexity. Increasing this degree from $d=2$ to $d=30$, the time required for parsing and solving the program using the \texttt{dpvar}, \texttt{pvar}, \texttt{syms} and \texttt{sdpvar} (YALMIP) implementations was determined. The results are displayed in Fig.~\ref{fig:polynomial_optimization}.

Solving the GLB problem with SOSTOOLS 3.04, the parsing complexity increases rapidly with the monomial degree, already exceeding a computation time of one hour for monomial degrees 10 (\texttt{syms}) or 12 (\texttt{pvar}). This rate of increase is substantially improved using SOSTOOLS 4.00, displaying a slope similar to that using YALMIP, though reducing computation time by a factor of around $10^2$. Moreover, the \texttt{dpvar} data structure is able to achieve a much more favorable solve-to-setup time ratio, with in general less than 20\% of the computation time spent on parsing.

\subsection{Robust Stability}\label{sec:subsec_Robust_Stability}

As a second example, we consider testing robust stability of a linear ODE 
\begin{align*}
\dot{x}(t)=A(p)x(t),
\end{align*}
with state $x(t)\in\R^n$ at any $t\geq 0$ and uncertain parameters $p\in G:=\{p\in\R^2\mid g(p)\geq 0\}$, where $g(p)=1-p_1^2-p_2^2$. Using a quadratic Lyapunov function $V(p,x)=x^T P(p)x$, we may determine stability of this system by testing for existence of a matrix-valued polynomial $P(p)$ such that $P(p)>0$ and $P(p)A(p)+A^T(p)P(p)\leq 0$ for any $p\in G$. Using the Psatz, we approach this as an SOS problem
\begin{align}\label{eq:robust_stability_SOS}
P-\epsilon I_n&\in \Sigma_s[p],	&
-Qg-PA-A^TP&\in \Sigma_s[p],
\end{align}
where $Q\in\Sigma_s[p]$, and we let $\epsilon=10^{-4}$.

In parsing this problem, we considered a polynomial matrix $A\in\R^{n\times n}[p]$ with all lower diagonal elements equal to $0.25p_1$, all upper diagonal elements equal to $-0.25p_2$, and all diagonal elements equal to $1$. 
The time required for parsing was computed for problem sizes up to $n=30$ and for each of the different implementations, using a variable $P$ of maximal degree $2d=4$. 
The results are displayed in Fig.~\ref{fig:robust_stability}.


The results again show that the \texttt{dpvar} implementation requires significantly less time to parse than the alternative implementations. This time also scales much more favorably using the \texttt{dpvar} data structure, in general offering an order $10^2$ reduction in computation time compared to all other implementations. In fact, even for $n=50$, the \texttt{dpvar} structure allowed the problem to be parsed in just 374 seconds, a threshold exceeded by YALMIP at $n=13$.



\subsection{Local Stability}\label{sec:subsec_Local_Stability}

As a final example, we test local stability of a chain of $n$ Van der Pol oscillators. In particular, we consider the system presented in~\cite{tacchi2020LocalStability}, given by $\dot{x}(t)=f(x)$,
where $x=(y,z)=(y_1,\hdots,y_n,z_1,\hdots,z_n)$ and
{\small
\begin{align*}
f_i(y,z)&=-2z_i,		\hspace*{2.85cm}	\forall i\in\{1,\hdots,n\}	\\
f_{n+j}(y,z)&=0.8y_j+10(1.2^2y_j^2-0.21)z_j+\epsilon_j z_{j+1}y_j,\\	&\hspace*{4.10cm}	\forall j\in\{1,\hdots,n-1\}	\\
f_{2n}(y,z)&=0.8y_n+10(1.2^2 y_n^2-0.21)z_n,
\end{align*}
}
where we let $\epsilon_j=-0.5$ for each $j$. We test stability inside a ball of radius $r=0.5$, so that $x\in \{x\in\R^{2n}\mid g(x)\geq 0\}$, where $g(x)=r^2-\|x\|^2$. To this end, we once again use a Lyapunov function $V\in\Sigma_s[x]$, imposing a Psatz condition
\begin{align}\label{eq:local_stability_SOS}
-[\nabla V(x)]^Tf(x)-s(x)g(x)\in\Sigma_s[x]
\end{align}
where $s\in\Sigma_s$.
Parsing this problem for increasing values of $n$, we once more determined the time required for parsing and solving the problem using the different implementations, using a function $V$ of degree $2d=4$. The results are presented in Fig.~\ref{fig:local_stability}. Note that the SDP solver ran out of memory for problems involving more than $2n=14$ independent variables, prohibiting further tests.

Solving the local stability problem with both the \texttt{pvar} and \texttt{syms} implementations, the required time to parse the SOS program almost consistently accounts for more than 80\% of the total computation time. This issue is resolved using the \texttt{dpvar} data structure, allowing SOSTOOLS 4.00 to parse the problem with an efficiency similar to that of YALMIP.

\section{Conclusion}

In this paper, we have introduced a new representation of polynomial variables, which is affine in the decision variables. We showed that, using this dpvar representation, computation time for polynomial operations such as addition, multiplication and differentiation remains relatively small, increasing favorably with the number of involved decision variables. Exploiting the MATLAB built-in sparsity structure, we also showed that the computational and memory overhead for storing and manipulating variables in the dpvar representation is minimal, allowing for efficient parsing of SOS programs. Incorporating this representation in SOSTOOLS 4.00, performance of this parser was drastically enhanced, requiring computation times similar to or even less than those using the batch parser YALMIP to parse common optimization problems.

\bibliographystyle{IEEEtran}
\bibliography{bibfile}

\newpage
\appendix

\subsection{Computational Complexity of Merging Monomial Bases}\label{apndx:merging_bases}

Representing polynomial variables using either the \texttt{pvar} or \texttt{dpvar} data formats, almost all binary operations require the monomial bases of the considered polynomial variables to be merged. For example, recall from Subsection~\ref{sec:subsec_addition} the problem of adding two (scalar) polynomial variables $s_1\in\R[x_1,\hdots,x_{p_1};\xi_1,\hdots,\xi_{q_1}]$ and $s_2\in\R[y_1,\hdots,y_{p_2};\eta_1,\hdots,\eta_{q_2}]$, written in the dpvar representation as
\begin{align*}
s_1(x;\xi)&=Z_1(\xi)^T C_1 Z_{d_1}(x),	\\
s_2(y;\eta)&=Z_1(\eta)^T C_2 Z_{d_2}(y).
\end{align*}
It is clear that the sum $s_3=s_1+s_2$ of these polynomials may be represented as
\begin{align*}
s_3(x,y;\xi,\eta)&=\bmat{Z_1(\xi)\\Z_1(\eta)}^T \bmat{C_1&0\\0&C_2} \bmat{Z_{d_1}(x)\\Z_{d_2}(y)}.
\end{align*}
To express this result in the dpvar representation, we have to define the variables $z,\chi$, monomial basis $\hat{Z}_{d_3}\in\R^{n_3}[z]$, and coefficients $C_3$ such that
\begin{align*}
s_3(z;\chi)&=Z_1(\chi)^T C_3 \hat{Z}_{d_3}(z)=\bmat{1\\\chi}^T C_3\hat{Z}_{d_3}(z).
\end{align*}
Here, merging the bases $Z_{d_1}(x)$ and $Z_{d_2}(y)$ into a single (incomplete) basis $\hat{Z}_{d_3}(z)$ of monomials of degree at most $d_3:=\max\{d_1,d_2\}$ in variables $z=\text{unique}(x;y)$ requires significant computational effort, often accounting for the greatest computational cost in performing operations like addition.

To get an estimate of the complexity associated with merging the bases, let $Z_{d_1}\in\R^{n_1}[x_1,\hdots,x_{p_1}]$ and $Z_{d_2}\in\R^{n_2}[y_1,\hdots,y_p]$ consist of respectively $n_1$ and $n_2$ monomials, in respectively $p_1$ and $p_2$ variables. The bases can then be represented as matrices $Z_{\text{M},d_1}\in\N^{n_1\times p_1}$ and $Z_{\text{M},d_2}\in\N^{n_2\times p_2}$ containing the degrees of each variable in each monomial, so that the full vector of monomials $\smallbmat{Z_{\text{M},d_1}(x)\\Z_{\text{M},d_2}(y)}$ can be represented by the matrix
\begin{align*}
\overbrace{\bmat{Z_{\text{M},d_1}&0\\0&Z_{\text{M},d_2}}}^{[\ x\ ,\ y\ ]}\in\N^{(n_1+n_2)\times(p_1+p_2)}.
\end{align*}
Conversion of this matrix into a degree matrix $\hat{Z}_{\text{M},d_3}\in\N^{n_3\times p_3}$ for the merged basis $\hat{Z}_{d_3}\in\R^{n_3}[z]$ is performed in 3 steps.

\subsubsection{Merging the variables}
First, a unique set of variables $z_1,\hdots,z_{p_3}$ is determined from $x_1,\hdots,x_{p_1}$ and $y_1,\hdots,y_{p_2}$. This can be done very efficiently using e.g. a quicksort algorithm to sort the variables, and discarding redundant appearances of each variable, requiring a cost of
\begin{align*}
 \mcl{O}\bl((p_1+p_2)\log(p_1+p_2)\br).
\end{align*}
In defining these variables $z$, we also obtain permutation matrices $P_1\in\N^{p_1\times p_3}$ and $P_2\in\N^{p_2\times p_3}$ such that
\begin{align*}
 \bmat{x_1\\\vdots\\x_{p_1}}&=P_1\bmat{z_1\\\vdots\\ z_{p_3}},	&
 &\text{and}	&
 \bmat{y_1\\\vdots\\y_{p_2}}&=P_2\bmat{z_1\\\vdots\\ z_{p_3}}.
\end{align*}
Using these permutation matrices, the full vector of monomials $\smallbmat{Z_{d_1}(x)\\Z_{d_2}(y)}$ may be equivalently represented by the degree matrix
\begin{align*}
\bmat{\hat{Z}_{d_1}\\\hat{Z}_{d_2}}
=\bmat{Z_{\text{M},d_1}P_1\\Z_{\text{M},d_2}P_2}
\in\N^{(n_1+n_2)\times p_3},
\end{align*}
describing the degrees of each monomial in terms of the new variables $z$.

\subsubsection{Sorting the monomials}
Next, the rows of $\smallbmat{\hat{Z}_{\text{M},d_1}\\\hat{Z}_{\text{M},d_2}}$ are ordered in lexicographical order. For this, a weight is assigned to each monomial, collected in a vector $\hat{\mbf{z}}\in\N^{n_1+n_2}$, computed as
\begin{align}\label{eq:apdx:mon_weight}
 \hat{\mbf{z}}=\bmat{[\hat{Z}_{\text{M},d_1}]_1&[\hat{Z}_{\text{M},d_1}]_2&\hdots [\hat{Z}_{\text{M},d_1}]_{p_3}\\
 [\hat{Z}_{\text{M},d_2}]_1&[\hat{Z}_{\text{M},d_2}]_2&\hdots [\hat{Z}_{\text{M},d_2}]_{p_3}}\bmat{(d_3+1)^{p_3}\\(d_3+1)^{(p_3-1)}\\\vdots\\(d_3+1)^1}.
\end{align}
Here, $[\hat{Z}_{\text{M},d_i}]_k\in\N^{n_i}$ denotes column $k$ of $\hat{Z}_{\text{M},d_i}\in\N^{n_i\times p_3}$, and $d_3:=\max\{d_1,d_2\}$ is the maximal degree of all monomials, so that $[\hat{Z}_{\text{M},d_i}]_{jk}<d_3+1$ for any $j\in\{1,\hdots,n_i\}$ and $k\in\{1,\hdots,p_3\}$. This ensures that $\hat{\mbf{z}}_j> \hat{\mbf{z}}_i\in\N $ for $i,j\in\{1,\hdots,n_1+n_2\}$ if and only if row $j$ of $\smallbmat{\hat{Z}_{\text{M},d_1}\\\hat{Z}_{\text{M},d_2}}$ is greater than row $i$ of this matrix in a lexicographical sense. The vector $\hat{\mbf{z}}$ is then sorted calling the MATLAB inherent function \texttt{sort}, applying the quicksort algorithm, invoking a complexity of
\begin{align*}
 \mcl{O}\bl((n_1+n_2)\log(n_1+n_2)\br).
\end{align*}
Sorting the monomials, we obtain a permutation matrix $P_{\text{sort}}\in\N^{(n_1+n_2)\times(n_1+n_2)}$ so that $\tilde{Z}_{\text{M},d_3}:=P_{\text{sort}}\smallbmat{\hat{Z}_{\text{M},d_1}\\\hat{Z}_{\text{M},d_2}}\in\N^{(n_1+n_2)\times p_3}$ contains the degrees of all monomials in lexicographical order.

\subsubsection{Discarding duplicate monomials}
Finally, a unique set of monomials can be obtained from the ordered set by comparing subsequent rows of the matrix $\tilde{Z}_{\text{M},d_3}$, retaining only the first of each pair $[\tilde{Z}_{\text{M},d_3}]_j=[\tilde{Z}_{\text{M},d_3}]_{j+1}$ of identical rows. Since the degrees are stored as a sparse matrix, only nonzero values need to be compared, resulting in a complexity
\begin{align*}
 \mcl{O}\bl(nnz(\tilde{Z}_{\text{M},d_3})\br)
 &=\mcl{O}\bl(nnz(\hat{Z}_{\text{M},d_1})+nnz(\hat{Z}_{\text{M},d_2})\br)	\\
 &=\mcl{O}\bl(nnz(Z_{\text{M},d_1})+nnz(Z_{\text{M},d_2})\br).
\end{align*}
We obtain a matrix $P_{\text{unique}}\in\N^{n_3\times (n_1+n_2)}$ such that
\begin{align*}
 \hat{Z}_{\text{M},d_3}:=P_{\text{unique}}\tilde{Z}_{\text{M},d_3}=P_{\text{unique}}P_{\text{sort}}\bmat{\hat{Z}_{\text{M},d_1}\\\hat{Z}_{\text{M},d_2}}\in\N^{n_3\times p_3}
\end{align*}
is a matrix of degrees associated to the unique combination of monomials in $Z_{d_1}(x)$ and $Z_{d_2}(y)$.

In performing these steps, it is clear that the sorting (Step 2) and subsequent comparing (Step 3) of the monomials $\bar{Z}_{d_3}\in\N^{(n_1+n_2)\times p_3}$ will require the greatest computational effort. We note here that, for $p_i$ variables and a maximal degree $d_i$, the total number $n_i$ of possible monomials is 
\begin{align*}
 n_i=\frac{(p_i+d_i)!}{p_i!d_i!}.
\end{align*}
Moreover, the number of nonzero elements in the degree matrix $Z_{\text{M},d_i}\in\N^{n_i\times p_i}$ associated to these monomials is given by
\begin{align*}
nnz(Z_{\text{M},d_i})&=\frac{(p_i+d_i)!-p_i\ [(p_i-1+d_i)!]}{(p_i-1)!\ d_i!}	\\
&=\left[p_i-\frac{p_i^2}{p_i+d_i}\right]n_i
\end{align*}
For sufficiently large values of $p_i$ and $d_i$, here,
\begin{align*}
	\left[p_i-\frac{p_i^2}{p_i+d_i}\right]\leq \log\left(\frac{(p_i+d_i)!}{p_i!d_i!}\right)=\log(n_i),
\end{align*}
and thus, in general, the complexity of sorting the monomials in $\smallbmat{\hat{Z}_{\text{M},d_1}\\\hat{Z}_{\text{M},d_2}}$ will be greater than that of merging duplicate monomials in the sorted $\tilde{Z}_{\text{M},d_3}$. We conclude that the complexity of merging the monomial bases $Z_{d_1}\in\R^{n_1}[x]$ and $Z_{d_2}\in\R^{n_2}[y]$ is roughly
\begin{align*}
\mcl{O}\bl((n_1+n_2)\log(n_1+n_2)\br).
\end{align*}
Here, $n_i:=\frac{(p_i+d_i)!}{p_i!d_i!}$, so that the cost of adding two polynomial variables increases rapidly with the number of independent variables $p_1$ and $p_2$. In this sense, the dpvar representation offers a significant advantage over the pvar representation, by not storing decision variables as independent variables, and thus maintaining relatively small values for $p_i$. 

It should be noted that the monomial sorting of $\smallbmat{\hat{Z}_{\text{M},d_1}\\\hat{Z}_{\text{M},d_2}}\in\N^{(n_1+n_2)\times p_3}$ described in Step 2, may require additional steps when considering large numbers of independent variables. In particular, for large values of $p_3$ and $d_3$, the weights $\hat{\mbf{z}}_j$ of each monomial, computed as in Equation~\eqref{eq:apdx:mon_weight}, may exceed the maximal numerical values MATLAB can (effectively) handle. Under these circumstances, sorting may have to be performed in stages, sorting only based on a subset of the columns of $\smallbmat{\hat{Z}_{\text{M},d_1}\\\hat{Z}_{\text{M},d_2}}$ at each stage. This will increase the complexity with a factor dependent on the number of stages in which the sorting has to be performed. This additional complexity is in general avoided when using the dpvar representation, as the number of variables and monomial degree in common SOS programs are usually sufficiently small. However, using the pvar representation, since the decision variables are included as independent variables in the monomial, the number of columns $p_3$ will be drastically increased, thus requiring further computational effort that can be avoided with the dpvar representation.

\subsection{A SOSTOOLS Implementation of Several Polynomial Optimization Problems}\label{apndx:optimization_tests}

\subsubsection{Greatest Lower Bound}

The greatest lower bound problem from Subsection~\ref{sec:subsec_Greatest_Lower_Bound} takes the form
\begin{align*}
\max_{\gamma}\quad \gamma,&	\\
\text{s.t.}\quad \gamma&\leq f(x) \quad
\forall x_1,x_2\in[-12,12],
\end{align*}
where $f(x)=x_1^4+x_2^4-2x_2x_1^3-3x_2^2x_1^2+150(x_1^2+x_2^2)$.
Defining,
\begin{align*}
g_1(x)&=12^2-x_1^2\geq 0,	&
g_2(x)&=12^2-x_2^2\geq 0,\\
g_3(x)&=2\cdot 12^2-(x_1^2+x_2^2)\geq 0,
\end{align*}
and invoking Putinar's Positivstellensatz (Psatz)~\cite{putinar1993Psatz} (Psatz), we enforce a single SOS constraint
\begin{align*}
F(x):=(f(x)-\gamma) - s_1g_1(x) - s_2g_2(x) - s_3g_3(x)\in\Sigma_s,
\end{align*}
with SOS variables $s_1,s_2,s_3\in\Sigma_s$. This SOS problem may be implemented in SOSTOOLS 4.00 by first initializing a program structure \texttt{sos} in the independent variables $x_1,x_2$ and decision variable $\gamma$, as
\begin{verbatim}
> pvar x1 x2
> dpvar gam
> sos = sosprogram([x1,x2],gam);
\end{verbatim}
Note here that the independent variables $x_1,x_2$ are implemented as \texttt{polynomial} (\texttt{pvar}) class objects, whereas the decision variable $\gamma$ is implemented as a \texttt{dpvar} class object. Next, SOS variables $s_i(x;C)=Z_{d}(x)C_i Z_{d}(x)$ for each $i\in\{1,2,3\}$ are initialized as,
\begin{verbatim}
> Zd = monomials([x1;x2],0:d)
> [sos,s1] = sossosvar(sos,Zd);   
> [sos,s2] = sossosvar(sos,Zd);  
> [sos,s3] = sossosvar(sos,Zd);  
\end{verbatim}
where now \texttt{Zd} will be a \texttt{polynomial} class object, representing a monomial vector $Z_{d}(x)$ of maximal degree $d$, and \texttt{si} will be \texttt{dpvar} class objects. Implementing the functions $f$ and $g_i$ as \texttt{polynomial} class objects $\texttt{f}$ and $\texttt{gi}$, the SOS constraint $F\in\Sigma_s$ is finally imposed as
\begin{verbatim}
> F = f-gam - s1*g1 - s2*g2 - s3*g3;
> sos = sosineq(sos,F);
\end{verbatim}
at which point the program can be solved by calling
\begin{verbatim}
> sos = sossolve(sos);
\end{verbatim}

\subsubsection{Robust Stability}

In Subsection~\ref{sec:subsec_Robust_Stability}, we consider a linear ODE
\begin{align*}
\dot{x}(t)=A(p)x(t),
\end{align*}
with state $x(t)\in\R^n$ at any $t\geq 0$ and uncertain parameters $p\in G:=\{p\in\R^2\mid g(p)\geq 0\}$, where $g(p)=1-p_1^2-p_2^2$. Robust stability is determined by testing for existence of a matrix-valued polynomial $P(p)$ such that $P(p)>0$ and $P(p)A(p)+A^T(p)P(p)\leq 0$ for any $p\in G$, enforced as an SOS problem
\begin{align*}
P-\epsilon I_n&\in \Sigma_s[p],	&
-Qg-PA-A^TP&\in \Sigma_s[p],
\end{align*}
where $Q\in\Sigma_s[p]$, and we let $\epsilon=10^{-4}$. 
In SOSTOOLS 4.00, after initializing an SOS program as
\begin{verbatim}
> pvar p1 p2
> sos = sosprogram([p1,p2]);
\end{verbatim}
the robust stability test may be implemented by first defining the positive definite polynomial variable $P\in\Sigma_s[p]$ in terms of monomials of degree 2 as
\begin{verbatim}
> Z=monomials([p1;p2],0:2)
> [sos,P]=sospolymatrixvar(sos,Z,[n n]);
> eps=1e-4;
> [sos]=sosmatrixineq(sos,P-eps*eye(n));
\end{verbatim}
where now \texttt{P} is a \texttt{dpvar} class object representing the SOS variable $P(p;C)=Z_2(p)^T C Z_2(p)$, and satisfying $P-\epsilon I\in\Sigma_s[p]$.
Next, defining \texttt{polynomial} class objects \texttt{A} and \texttt{g} to represent the functions $A(p)$ and $g(p)$ respectively, negativity of the derivative is enforced as
\begin{verbatim}
> [sos,Q]=sospolymatrixvar(sos,Z,[n n]);
> [sos]=sosmatrixineq(sos,Q);
> [sos]=sosmatrixineq(sos,-Q*g-A’*P-P*A);
\end{verbatim}
at which point the program can be solved by calling
\begin{verbatim}
> sos = sossolve(sos);
\end{verbatim}

\subsubsection{Local Stability}

In Subsection~\ref{sec:subsec_Local_Stability}, we consider a system presented in~\cite{tacchi2020LocalStability}, given by $\dot{x}(t)=f(x)$,
where $x=(y,z)=(y_1,\hdots,y_n,z_1,\hdots,z_n)$ and
\begin{align*}
f_i(y,z)&=-2z_i,		\hspace*{3cm}	\forall i\in\{1,\hdots,n\}	\\
f_{n+j}(y,z)&=0.8y_j+10(1.2^2y_j^2-0.21)z_j+\epsilon_j z_{j+1}y_j,\\	&\hspace*{3.75cm}	\forall j\in\{1,\hdots,n-1\}	\\
f_{2n}(y,z)&=0.8y_n+10(1.2^2 y_n2-0.21)z_n,
\end{align*}
where we let $\epsilon_j=-0.5$ for each $j$. Local stability of this system is tested inside a ball of radius $r=0.5$, so that $x\in G:=\{x\in\R^{2n}\mid g(x)\geq 0\}$, where $g(x)=r^2-\|x\|^2$. To this end, a Lyapunov function $V\in\Sigma_s[x]$ is sought, imposing an SOS constraint
\begin{align*}
-[\nabla V(x)]^Tf(x)-s(x)g(x)\in\Sigma_s[x]
\end{align*}
where $s\in\Sigma_s$. This SOS problem may be implemented as a program structure \texttt{sos} in SOSTOOLS, initialized as
\begin{verbatim}
> pvar y1 ... yn;
> pvar z1 ... zn;
> sos = sosprogram([y1,...,zn]);
\end{verbatim}
Next, we construct a variable $V(x;C)=Z_2(x)^T C Z_2(x)$,
\begin{verbatim}
> Z = monomials([y1,...,zn],0:2);
> [sos,V] = sossosvar(sos,Z);
\end{verbatim}
defining a \texttt{dpvar} class object \texttt{V} representing the Lyapunov function. Defining \texttt{polynomial} class objects \texttt{f} and \texttt{g} to represent the desired functions $f(x)$ and $g(x)$, the derivative of the Lyapunov function is finally enforced to be negative in the desired domain
\begin{verbatim}
> Vd = jacobian(V,[y1,...,zn])*f;
> [sos,s] = sossosvar(sos,Z);
> [sos] = sosineq(sos,-Vd-s*g);
\end{verbatim}
at which point the program can be solved by calling
\begin{verbatim}
> sos = sossolve(sos);
\end{verbatim}

\subsection{The \texttt{sosquadvar} Function}\label{apndx:sosquadvar}

In addition to incorporating the \texttt{dpvar} data structure, SOSTOOLS 4.00 also introduces the \texttt{sosquadvar} function, for efficient implementation of general polynomial decision variables. In its simplest form, \texttt{sosquadvar} takes as input a SOSTOOLS program structure \texttt{sos}, and two monomial vectors $Z_{d_1}\in\R^{k_1}[x]$ and $Z_{d_2}\in\R^{k_2}[y]$, implemented as \texttt{polynomial} (\texttt{pvar}) class objects \texttt{Z1} and \texttt{Z2}. Calling
\begin{verbatim}
> [sos,P] = sosquadvar(sos,Z1,Z2);
\end{verbatim}
a \texttt{dpvar} class object \texttt{P} is returned, representing a polynomial variable $P(x,y;Q)=Z_{d_1}(x)^T Q Z_{d_2}(y)$, for decision variables $Q\in\R^{k_1\times k_2}$. The decision variables are also added to the output program structure \texttt{sos}. Using the \texttt{sosquadvar} function, monomial vectors $Z_{d_1}=1$ or $Z_{d_2}=1$ may also be specified, allowing e.g. linear polynomial variables $P(y;q)=q^T Z_{d_2}(y)$ to be added to the program. Moreover, optional matrix dimensions \texttt{m} and \texttt{n} may also be passed to the function as
\begin{verbatim}
> [sos,P] = sosquadvar(sos,Z1,Z2,m,n);
\end{verbatim}
producing a \texttt{dpvar} object \texttt{P} associated to the $m\times n$ matrix-valued variable 
\[
P(x,y;Q)=(I_{m}\otimes Z_{d_1}(x))^T Q (I_{n}\otimes Z_{d_2}(y)),
\]
where now $Q\in\R^{m k_1\times nk_2}$.

In addition to the dimensions of the variable, positivity properties of the variable can be specified when calling \texttt{sosquadvar}. In particular, the function allows a sixth (optional) input to be passed, taking one of two values:
\begin{enumerate}
	\item `\texttt{sym}', requiring the decision variable $Q\in\R^{m k_1\times nk_2}$  to be symmetric, or
	\item `\texttt{pos}', requiring the decision variable $Q\in\R^{m k_1\times nk_2}$ to be (symmetric) positive semi-definite.
\end{enumerate}
Naturally, both of these options only make sense if the matrix $Q$ is square, allowing these options to be specified only if $m=n$ and $k_1=k_2$. Using the \texttt{pos} input, an SOS variable $S(x;Q)=(I_{m}\otimes Z_{d_1}(x))^T Q (I_{m}\otimes Z_{d_1}(x))$ with $Q\geq 0$ can be added to the program by calling
\begin{verbatim}
> [sos,S]=sosquadvar(sos,Z1,Z1,m,m,'pos');
\end{verbatim}
In calling the function with this \texttt{pos} input, the constraint $Q\geq 0$ on the decision variables of $S(x;Q)$ will be added to the program structure \texttt{sos}. Note, however, that unless the left and right monomial vectors are identical, the resulting variable $S(x,y;Q)$ need not be an SOS variable. 

As a final functionality, \texttt{sosquadvar} allows variables to be specified for which positivity is coupled between multiple polynomial variables. Specifically, consider two sets $\{Z_{d_{1,1}},\hdots,Z_{d_{1,r}}\}$ and $\{Z_{d_{2,1}},\hdots,Z_{d_{2,p}}\}$ of respectively $r\in\N$ and $p\in\N$ monomial vectors, where $Z_{d_{1,i}}\in\R^{k_{1,i}}[x_i]$ and $Z_{d_{2,j}}\in\R^{k_{1,j}}[x_j]$ for each $i\in\{1,\hdots,r\}$ and $j\in\{1,\hdots,p\}$. For each pair of monomials $(Z_{d_{1,i}},Z_{d_{2,j}})$, \texttt{sosquadvar} can be used to construct a polynomial variable 
\begin{align*}
 P_{i,j}(x_i,y_j;Q_{i,j})&=(I_{m_i}\otimes Z_{d_{1,i}}(x_i))^T Q_{i,j} (I_{n_j}\otimes Z_{d_{2,j}}(y_j))	\\
 &\hspace*{2.0cm} \in\R^{m_i\times n_j}[x_i,y_j;Q_{i,j}],
\end{align*}
parameterized by decision variables $Q_{i,j}\in\R^{m_i k_{1,i}\times n_i k_{2,j}}$. Defining such variables for each pair $(i,j)$ separately, however, positivity of the matrices $Q_{i,j}$ is not necessary or sufficient for positivity of the composite matrix
\begin{align}\label{eq:quadvar_Qbigmat}
 Q=\bmat{Q_{1,1}&\hdots &Q_{1,p}\\
 	\vdots &\ddots &\vdots\\
 	Q_{r,1}&\hdots & Q_{r,p}}\in\R^{\sum_{i=1}^{r}m_i k_{1,i}\times \sum_{j=1}^{p}m_j k_{2,j}}
\end{align}
as a whole. Instead, to construct the polynomials $P_{i,j}\in\R^{m_i\times n_j}[x_i,y_j;Q_{i,j}]$ while enforcing $Q\geq 0$, \texttt{sosquadvar} can be called with MATLAB cell structures \texttt{Z1=\{Z11,...,Z1r\}} and \texttt{Z2=\{Z21,...,Z2p\}}, where \texttt{Z1i} and \texttt{Z1j} are \texttt{polynomial} class objects defining the desired monomial vectors $Z_{d_{1,i}}$ and $Z_{d_{2,j}}$. Using vectors \texttt{m=[m1,...,mr]} and \texttt{n=[n1,...,np]} to specify the matrix dimensions, \texttt{sosquadvar} can be called as before,
\begin{verbatim}
> [sos,P]=sosquadvar(sos,Z1,Z2,m,n,'pos');
\end{verbatim}
producing an $r\times p$ cell structure \texttt{P}, where each element \texttt{P\{i,j\}} is a \texttt{dpvar} class object representing the polynomial variable $P_{i,j}\in\R^{m_i\times n_j}[x_i,y_j;Q_{i,j}]$, and where the matrix $Q$ as in Eqn.~\eqref{eq:quadvar_Qbigmat} is required to satisfy $Q\geq 0$. Calling \texttt{sosquadvar} with cell inputs, the \texttt{pos} and \texttt{sym} options can only be used if $r=p$, and $m_i=n_i$ and $k_{1,i}=k_{2,i}$ for each $i\in\{1,\hdots,r\}$. If for each $i$ further $Z_{d_{1,i}}=Z_{d_{2,i}}$, and the \texttt{pos} option is specified, the composite variable $P\in\R^{\sum_{i=1}^{r}m_i\times \sum_{i=1}^{r}n_i}[x;Q]$ will be an SOS variable, though the individual functions $P_{i,j}(x_i,x_j;Q_{i,j})$ (for $i\neq j$) will generally not be.

Through the \texttt{sosquadvar} function, SOSTOOLS 4.00 allows straightforward implementation of a wide class of polynomial variables, substantially expanding the scope of variables that could be specified in SOSTOOLS 3.04. Constructing these variables directly as \texttt{dpvar} objects, \texttt{sosquadvar} also increases efficiency compared to the functions \texttt{sossosvar}, \texttt{sospolyvar}, \texttt{sosposmatrvar}, etc., used for constructing different types of polynomial variables in SOSTOOLS 3.04. Accordingly, each of these functions has been updated to outsource computations to \texttt{sosquadvar} where possible, enhancing efficiency and transparency in the parsing of SOS programs.

\end{document}